\documentclass[article, floatsintext]{elsarticle}
\journal{}

\usepackage[english]{babel}
\usepackage[utf8]{inputenc}
\usepackage[displaymath, mathlines]{lineno}

\usepackage{amsfonts}
\usepackage{amsmath}
\usepackage{amssymb}
\usepackage{amsthm}
\usepackage{bm}
\usepackage{bbm}
\usepackage{booktabs}
\usepackage[labelsep=period]{caption}
\usepackage{enumitem}
\usepackage{etoolbox}
\usepackage{fancyhdr}
\usepackage{float}
\usepackage{geometry}
\usepackage{graphicx}
\usepackage{lastpage}
\usepackage{listings}
\usepackage{longtable}
\usepackage{subcaption}
\usepackage{units}
\usepackage[svgnames]{xcolor}

\newtheorem{theorem}{Theorem}[section]

\newtheorem{definition}[theorem]{Definition}

\definecolor{Yellow}{RGB}{255, 192, 0}
\definecolor{Orange}{RGB}{237, 125, 49}
\definecolor{Purple}{RGB}{112, 48, 160}
\definecolor{Raspberry}{RGB}{212, 52, 117}
\definecolor{Blue}{RGB}{0,161,218}
\definecolor{Green}{RGB}{112, 173, 71}
\definecolor{Teal}{RGB}{5,145,130}

\newcommand{\Caratheodory}{{Carath\'eodory}}

\newcommand{\realR}{\mathbb{R}}

\DeclareMathOperator{\diag}{diag}

\DeclareMathOperator{\conv}{conv}

\newcommand{\bhat}[1]{\hat{\bm{#1}}}
\newcommand{\dint}[1]{~\text{d}#1}

\newcommand{\domain}{\Omega}
\newcommand{\CartDomain}{\hat{\Omega}}
\newcommand{\embeddedRegions}{\Omega_E}

\newcommand{\Nspace}{{N_d}}

\newcommand{\totalDegree}{\mathbb{P}}
\newcommand{\tensorProduct}{\mathcal{Q}}

\newcommand{\solnSpace}{\mathcal{V}}
\newcommand{\basisp}{\varphi}

\newcommand{\mesh}{\mathcal{T}}
\newcommand{\element}{D}
\newcommand{\refElement}{\hat{D}}

\newcommand{\entropyDomain}{H}
\newcommand{\entropy}{U}
\newcommand{\entropyFlux}{F}
\newcommand{\entropyPotential}{\psi}

\newcommand{\Q}{\bm{Q}}
\newcommand{\D}{\bm{D}}

\newcommand{\M}{\bm{M}}
\newcommand{\B}{\bm{B}}
\newcommand{\bP}{\bm{P}}
\newcommand{\E}{\bm{E}}
\newcommand{\ones}{\mathbbm{1}}

\newcommand{\f}{\bm{f}}
\newcommand{\F}{\bm{F}}

\newcommand{\fEC}{\bm{f}_{EC}}
\newcommand{\fES}{\bm{f}_{ES}}

\newcommand{\uAux}{\tilde{\bm{u}}}

\newcommand{\bx}{\bm{x}}
\newcommand{\bX}{\bm{X}}

\newcommand{\bu}{\bm{u}}

\newcommand{\bv}{\bm{v}}
\newcommand{\bV}{\bm{V}}
\newcommand{\bw}{\bm{w}}
\newcommand{\bW}{\bm{W}}

\newcommand{\bn}{\bm{n}}

\newcommand{\pderiv}[2]{\dfrac{\partial #1}{\partial #2}}
\newcommand{\deriv}[2]{\dfrac{\text{d}#1}{\text{d}#2}}

\newcommand{\pd}[3]{\frac{\partial^{#3} #1}{\partial#2^{#3}}}

\newcommand{\jump}[1] {\ensuremath{\LRs{\![#1]\!}}}
\newcommand{\avg}[1] {\ensuremath{\LRc{\!\{#1\}\!}}}

\newcommand{\LRp}[1]{\left( #1 \right)}
\newcommand{\LRs}[1]{\left[ #1 \right]}

\newcommand{\LRb}[1]{\left| #1 \right|}
\newcommand{\LRc}[1]{\left\{ #1 \right\}}

\begin{document} 

\begin{frontmatter}


\cortext[cor1]{Corresponding author(s)}

\author[1]{Christina G. Taylor\corref{cor1}}
\ead{christina.taylor@austin.utexas.edu}

\author[2]{Jesse Chan}
\ead{jesse.chan@rice.edu}

\affiliation[1]{organization={Oden Institute for Computational Science and Engineering, University of Texas at Austin}}
	
\affiliation[2]{organization={Department of Computational Applied Mathematics and Operations Research, Rice University}}


\title{An Entropy Stable High-Order Discontinuous Galerkin Method on Cut Meshes}
\date{Fall 2024}

\begin{keyword}
Entropy stable \sep Discontinuous Galerkin \sep Summation-by-parts \sep Cut meshes \sep Carath\'eodory pruning \sep State redistribution
\end{keyword}

\begin{abstract}
High-order entropy stable summation-by-parts (SBP) schemes are a class of robust and accurate numerical methods for hyperbolic conservation laws that are numerically stable at arbitrary order without the need for artificial stabilization. While SBP schemes are well-established on simplicial and tensor-product elements, they have not been extended to cut meshes. Cut meshes provide a convenient and efficient means of mesh generation for domains with embedded boundaries but can be difficult to use due to their arbitrarily shaped cut elements. Using the skew-hybridized SBP formulation of Chan \cite{chan-SkewSymmES}, we present a high-order accurate, entropy stable scheme for hyperbolic conservation laws on cut meshes. The formulation requires positive/non-negative weight quadrature rules on cut elements, which we construct via explicit parameterizations, subtriangulations, and Carath\'eodory pruning. We numerically verify the accuracy and stability of our method using the shallow water and compressible Euler equations and note promising results for the use of state redistribution with entropy stable methods.
\end{abstract}

\end{frontmatter}

\section{Introduction} 
We are interested in numerically solving systems of hyperbolic conservation laws taking the form
\begin{equation}
    \pderiv{\bu}{t} + \sum_{d=1}^{\Nspace} \pderiv{\f_d(\bu)}{x_d} = \bm{0}, \quad \bx \in \domain \subset \realR^\Nspace\label{eq:hyperbolicConsLaw} 
\end{equation}
with appropriate boundary and initial conditions and where the solution $\bu:\realR^n \to \entropyDomain$ for some convex set $\entropyDomain \subset \realR^n$. The functions $\f_d:\entropyDomain \to \realR^n, d=1,..., \Nspace$ are the fluxes and can be linear or nonlinear. A number of important physical systems fall into this framework, including the shallow water, compressible Euler, and ideal magento-hydrodynamics equations.

\begin{figure}[H]
    \centering
    \includegraphics[width=0.5\linewidth]{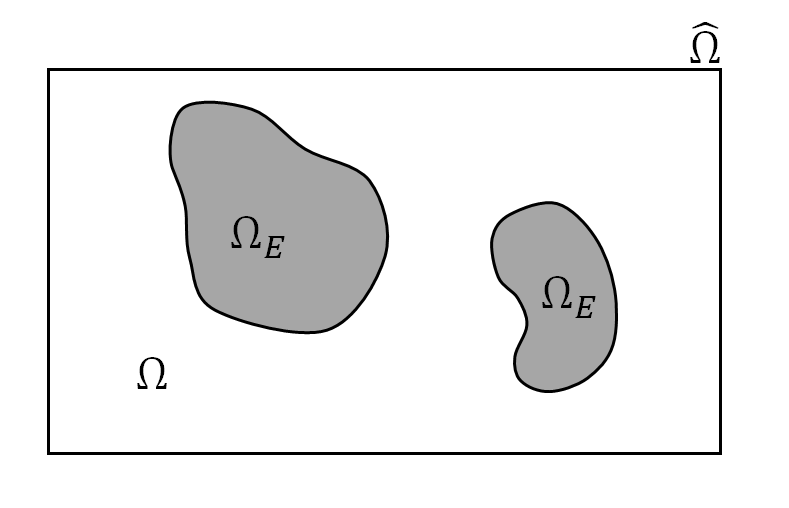}
    \caption{An example of a domain with embedded boundaries.}\label{fig:domainDef}
\end{figure}

We are interested in solving \eqref{eq:hyperbolicConsLaw} on domains $\domain \subset \realR^\Nspace$ featuring embedded boundaries. For this work we restrict ourselves to 2D domains, i.e., $\Nspace = 2$. The ability to represent embedded boundaries is fundamental to a number of fluid dynamics problems, which often consider flow past a object. Domains with such boundaries can be defined mathematically as a Cartesian embedding domain, $\CartDomain$, from which embedded regions, $\embeddedRegions$, are removed, i.e.
\begin{equation}
    \domain = \CartDomain - \embeddedRegions \label{eq:domainDef}.
\end{equation}
An example of such a domain is given in Figure \ref{fig:domainDef}. 

The main difficulty of this work is in ensuring entropy stability and high-order accuracy on a cut mesh. Hyperbolic conservation laws present a number of challenges to numerical approximation. First and foremost, hyperbolic conservation laws can develop shocks, i.e., discontinuities, in their solution even when starting from a smooth initial condition \cite{dafermos-HyperbolicConservationLawsBook}. Additionally, weak solutions of hyperbolic conservation laws need not be unique \cite{dafermos-HyperbolicConservationLawsBook} and some may not in fact be physical. 

Physically relevant weak solutions can be selected by augmenting \eqref{eq:hyperbolicConsLaw} with additional requirements on a solution $\bu$. One such set of requirements comes in the form of entropy stability. Despite their ability to develop discontinuities and sharp gradients (which can induce numerical instability), hyperbolic conservation laws are well-behaved in another sense: many satisfy a statement of entropy conservation or stability, a nonlinear analog of energy conservation/stability. Indeed, entropy conservation and stability can be seen as an extension of energy conservation/stability for linear systems \cite{gustafsson-L2-stability} to nonlinear systems. 

A statement of entropy conservation can easily derived from Equation \eqref{eq:hyperbolicConsLaw} given a scalar-valued, convex \textit{entropy function}, $\entropy:\entropyDomain \to \realR$ and \textit{entropy fluxes} $\entropyFlux_d:\entropyDomain \to \realR, d= 1,2$ satisfying
\begin{equation}
    \deriv{\entropyFlux_d}{\bu} = \deriv{\entropy}{\bu} \deriv{\f_d}{\bu}.
\end{equation}
Note here we use the convention that the derivative/partial derivatives of a function $g: \realR^n \to \realR^m$ with respect to its argument is in $\realR^{m \times n}$; this implies the derivative/gradient a scalar with respect to a column vector is a row vector and the derivative of a column vector with respect to a scalar is a column vector.

While the entropy function is an arbitrary convex function and thus (in general) non-unique, it is typically chosen to be a quantity with physical relevance, such as the total energy for the shallow water equations \cite{wintermeyer-SWE-flux} or the kinetic energy for Burgers equation. Stronger conditions on $\entropy$, such as uniform convexity or strict convexity, may yield other properties for a scheme or additional conditions on $\entropy$. 

In the case of uniformly convex $\entropy$, a priori $L_2$ bounds on entropy stable solutions are sometimes possible \cite{harten-ES-gudonov}. Strictly convex entropy functions are often of interest as the mapping between conservative variables $\bu$ and entropy variables $\bv$ is then one-to-one and can serve as a change of variables \cite{chen-ESFluxes}. For strictly convex $\entropy$, the entropy-entropy flux pair must symmetrize \eqref{eq:hyperbolicConsLaw} \cite{godlewski-hyperbolicMethodsBook, chen-ESFluxes}; this additional requirement can be enough to guarantee or suggest a unique entropy function for some systems, such as for the compressible Navier-Stokes equations with heat conduction \cite{hughes-NavierStokesFlux}. 

Historically, (and still often for the compressible Euler and Navier-Stokes equations \cite{harten-ESChainrule}) the entropy function was often related to the negative of physical, thermodynamic entropy, $S$, which established the convention of denoting the entropy function as $S(\bu)$. Here we follow the trend of more recent works in denoting the entropy function as $\entropy(\bu)$ to emphasize the arbitrary nature of the function, which need only be convex and have corresponding entropy fluxes.

To derive the statement of entropy conservation from \eqref{eq:hyperbolicConsLaw}, we define the entropy variables $\bv:\entropyDomain \to \realR^n$ as
\begin{equation}
    \bv(\bu) = \LRp{\deriv{\entropy}{\bu}}^T
\end{equation}
where we use the transpose to remain consistent with the notation used in the literature. By the chain rule and from the definition of the entropy function, entropy fluxes, and entropy variables we have that:
\begin{equation}
    \bv^T \pderiv{\bu}{t} = \deriv{\entropy}{\bu}{} \pderiv{\bu}{t}  = \pderiv{\entropy}{t}\label{eq:entropyRes-intro}
\end{equation}
and
\begin{equation}
    \bv^T \pderiv{\f_d}{x_d}{} = \deriv{\entropy}{\bu} \LRp{\deriv{\f_d}{\bu}\pderiv{\bu}{x_d}{}} = \LRp{\deriv{\entropy}{\bu} \deriv{\f_d}{\bu}}\pderiv{\bu}{x_d}{} = \deriv{\entropyFlux_d}{\bu}{}\pderiv{\bu}{x_d}{} =
    \pderiv{\entropyFlux_d}{x_d}{}.\label{eq:v-entropyFlux}
\end{equation}
Left multiplying Equation \eqref{eq:hyperbolicConsLaw} by $\bv^T$ and applying Equations \eqref{eq:entropyRes-intro} and  \eqref{eq:v-entropyFlux} yields a conservation law for the entropy,
\begin{equation}
    \pderiv{\entropy(\bu)}{t} + \sum_{d=1}^2 \pderiv{\entropyFlux_d(\bu)}{x_d} = 0,\label{eq:entropyConservation}
\end{equation}
which holds weakly for smooth $\bu$. In the presence of shocks, it can be shown via vanishing viscosity arguments that entropy is dissipated at shocks \cite{chen-ESFluxes}, yielding the entropy inequality:
\begin{equation}
    \pderiv{\entropy(\bu)}{t} + \sum_{d=1}^2 \pderiv{\entropyFlux_d(\bu)}{x_d} \leq 0\label{eq:entropyInequality}
\end{equation}
which is again to be interpreted in a weak sense. Inequality \eqref{eq:entropyInequality} in particular testifies of the stability of the solution $\bu$ with respect to a given entropy-entropy flux pair: it states that entropy can only decrease in time.

Numerical solutions to \eqref{eq:hyperbolicConsLaw} that are able to satisfy an appropriate discrete or semi-discrete analog of the entropy conservation law \eqref{eq:entropyConservation} are said to be \textit{entropy conservative}. Similarly, those satisfying analogs of the entropy inequality \eqref{eq:entropyInequality} are said to be \textit{entropy stable}. Importantly, entropy stability in numerical methods helps to prevent numerical instability. 

One of the first classes of entropy stable methods are based on the use of entropy stable numerical fluxes. The origins of flux-based entropy stable methods has its roots in first-order finite volume (FV) methods \cite{tadmor-fluxes, tadmor-EStheory, tadmor-ES-FV-FD-schemes, ray-ES-FV}. That original context was expanded into several finite volume methods and, via ``flux-differencing", discontinuous Galerkin (DG) methods. Some of these flux-based entropy stable methods can also be augmented by other properties/procedures such as essentially non-oscillatory (ENO) methods \cite{fjordholm-TeCNO} and kinetic energy preservation \cite{chandrashekar-KEP-ES-FV}. Fundamental to such methods are entropy conservative/stable two-point numerical fluxes as first defined by Tadmor in \cite{tadmor-fluxes}. 

The definition for entropy conservative and entropy stable two-point fluxes in the sense of Tadmor \cite{tadmor-fluxes} are given in Definitions \ref{def:ec-flux} and \ref{def:es-flux}. Necessary for these definitions are the \textit{entropy potentials}, $\entropyPotential_d:\entropyDomain \to \realR,$ which are defined as
\begin{equation}
    \entropyPotential_d(\bu) = \bv(\bu)^T\f_d(\bu) - F_d(\bu), \quad d = 1, 2.
\end{equation}
\begin{definition}
    In a given spatial dimension $\bx_d$, a two-point numerical flux $\fEC^d:\entropyDomain \times \entropyDomain \to \realR^n$ is said to be entropy conservative in the sense of Tadmor if for any $\bu_L, \bu_R \in \entropyDomain$ it satisfies:
    \begin{align}
        \fEC^d(\bu, \bu) &= \f_d(\bu) & \text{(consistency)} \\
        \fEC^d(\bu_L, \bu_R) &= \fEC^d(\bu_R, \bu_L) & \text{(symmetry)} \\
        \LRp{\bv_R - \bv_L}^T\fEC^d(\bu_L, \bu_R) &= \entropyPotential_d(\bu_R) - \entropyPotential_d(\bu_L) & \text{(entropy conservation)}
    \end{align}
    where $\bv_L = \bv(\bu_L), \bv_R = \bv(\bu_R)$.
\end{definition}\label{def:ec-flux}
\begin{definition}
    In a given spatial dimension $\bx_d$, a two-point numerical flux $\fES^d:\entropyDomain \times \entropyDomain \to \realR^n$ is said to be entropy stable in the sense of Tadmor if for any $\bu_L, \bu_R \in \entropyDomain$ it satisfies:
    \begin{align}
        \fES^d(\bu, \bu) &= \f_d(\bu) & \text{(consistency)} \\
        \LRp{\bv_R - \bv_L}^T\fES^d(\bu_L, \bu_R) &\leq \entropyPotential_d(\bu_R) - \entropyPotential_d(\bu_L) & \text{(entropy stability)}
    \end{align}
    where $\bv_L = \bv(\bu_L), \bv_R = \bv(\bu_R)$.
\end{definition}\label{def:es-flux}

Note in the case of entropy stable fluxes the symmetry requirement has been dropped and the entropy condition relaxed to an inequality. When performing calculations on the boundary of an element, information from the outward normal vector of the element is also needed to compute the final flux value in a given direction. When using entropy conservative fluxes, no changes are needed due to their symmetry property; the flux can simply be multiplied by the appropriate component of the outward normal vector. However, since entropy stable fluxes lack symmetry the classical two-point entropy stable flux $\fES^d(\bu_L, \bu_R)$ is often augmented to a three-argument flux $\fES^d(\bu_L, \bu_R, \bn_d)$ where $\bn_d$ is the $d^{th}$ component of the outward normal vector $\bn$. The inclusion of the normal vector induces a pseudo ``anti-symmetry" property
\begin{equation}
    \fES^d(\bu_L, \bu_R, \bn_d) = -\fES^d(\bu_R, \bu_L, -\bn_d)
\end{equation}
which is vital for the proof of global entropy stability of such schemes.

In practice, entropy stable fluxes can be constructed as a combination of a consistent, symmetric (but not necessarily entropy conservative) two-point base flux $\f_{CS}^d:\entropyDomain \times \entropyDomain \to \realR^n$ combined with a two-point dissipation term $\bm{d}_{\entropy}^d:\entropyDomain \times \entropyDomain \to \realR^n$:
\begin{equation}
    \fES^d(\bu_L, \bu_R, \bn_d) = \bn_d \f_{CS}(\bu_L, \bu_R) + \bm{d}_{\entropy}^d(\bu_L, \bu_R)\label{eq:ES-ds-flux}
\end{equation}
An example of such an entropy-stable flux is the Lax-Friedrichs flux
\begin{equation}
    \fES^d(\bu_L, \bu_R, \bn_d) = \bn_d\frac{1}{2}\LRp{\bu_L + \bu_R} - \frac{\lambda}{2}\LRp{\bu_R - \bu_L}
\end{equation}
which can also be expressed in shorthand notation as
\begin{equation}
    \fES^d(\bu_L, \bu_R) = \bn_d\avg{\bu} - \frac{\lambda}{2}\jump{\bu}
\end{equation}
where $\avg{\bu}$ is the arithmetic average of $\bu_L$ and $\bu_R$, $\jump{\bu} = \bu_R - \bu_L$ is the ``jump" of $\bu$, and $\lambda > 0$ is the maximum wavespeed or a suitable approximation thereof. 

The Lax-Friedrichs flux is a special case of a the HLL flux (named for the HLLE/HLLC Riemann solvers which are themselves named for Harten, Lax, and van Leer \cite{harten-HLL-solver}) and is entropy stable for suitable choices of $\lambda$ \cite{chen-ESFluxes, tan-laxWendroff, guermond-lambda-euler} despite its use of the average of the left and right states as the base flux. 

\subsection{High-Order Summation-By-Parts Operators}
While low-order methods can be extremely robust, their accuracy and efficiency per degree of freedom (DOF) can leave much to be desired. Low-order methods have higher numerical dispersion error than high-order equivalents \cite{ainsworth-highOrder-dispersionError} and are less efficient per degree of freedom \cite{wang-highVsLowOrder-HOinstability}. As a result, high-order methods are increasingly preferred for fluid applications that are especially sensitive to diffusion, such as unsteady flows featuring vortices and turbulence \cite{visbal-highOrder-unsteadyFlows}, and the simulation over long distances or times \cite{wang-highVsLowOrder-HOinstability}. However, the low-diffusion error of high-order methods is a double-edged sword: high-order methods tend to be unstable in the presence of shocks and discontinuities \cite{wang-highVsLowOrder-HOinstability}. 

In the context of hyperbolic conservation laws, instability can be attributed in part to a method's failure to satisfy the entropy inequality. Entropy stability can be imposed on a scheme externally by adding routines/operators that dissipate entropy to an otherwise entropy unstable scheme. Common approaches for externally enforcing entropy stability include artificial viscosity and regularization \cite{tominec-blastWave, nazarov-blastWave-original}, and filtering and/or limiting procedures \cite{zhang-maximumPrinciple, zhang-positivityPreservation}. However, these approaches can result in the loss of high-order accuracy \cite{wang-highVsLowOrder-HOinstability} and the need to tune parameters can result in a loss of robustness from one simulation to the next: too little dissipation and the scheme is unstable, too much and the solution is dominated by dissipation error.

A major difficulty in designing an inherently entropy stable high-order method is the use of the chain rule in space (as shown in \eqref{eq:v-entropyFlux}) in the derivation of the entropy statement as it has no semi-discrete equivalent. Summation-by-parts (SBP) operators paired with flux-differencing provides a means to side-step the need for the chain rule in the derivation of the semi-discrete entropy statement, a proof of which is presented in \cite{chan-hybridizedSBP}. SBP operators are so named as they are a semi-discrete analog of integration-by-parts and as such allow information in a volume to be replaced by information on the volume's boundary. 

SBP operators have been established for a number of finite volume \cite{tadmor-fluxes, tadmor-EStheory, fjordholm-TeCNO, chandrashekar-KEP-ES-FV, tadmor-ES-FV-FD-schemes, ray-ES-FV} and DG schemes \cite{fisher-ES-FD, carpenter-ES-spectralCollocation, gassner-splitFormDG, gassner-BR1}. SBP operators can also be constructed for finite difference schemes \cite{crean-SBP-generalCurvedElements, chen-ESFluxes}. SBP schemes are high-order accurate and can be used to create inherently entropy conservative/stable schemes.

Here, we use an entropy stable high-order SBP DG scheme, which pairs SBP operators with a flux-differencing DG formulation and entropy conservative/stable fluxes. In the DG setting, on a given element of the mesh classical diagonal-norm SBP schemes require the mass matrix $\M$ to be diagonal and positive-definite and SBP operators $\Q_d$ (which act as differentiation operators on volume of the element) and boundary integration matrices $\B_d$ satisfy the SBP property:
\begin{equation}
    \Q_d + \Q_d^T = \B_d, \quad d = 1,2.
\end{equation}

\subsection{Constructing Summation-By-Parts Operators on Cut Meshes}
While entropy stable high-order SBP methods are well-established, their treatment is typically limited to simplicial and tensor product elements; in 2D spatial domains these correspond to triangular and quadrilateral elements, respectively. Here, we are interested in domains with embedded boundaries. Embedded boundaries are most commonly handled via unstructured, fitted triangular meshes but such meshes can require either $h$ adaptivity or curvilinear elements to resolve curved boundaries well without loss of high-order accuracy. 

\begin{figure}[H]
    \centering
    \includegraphics[width=0.5\linewidth]{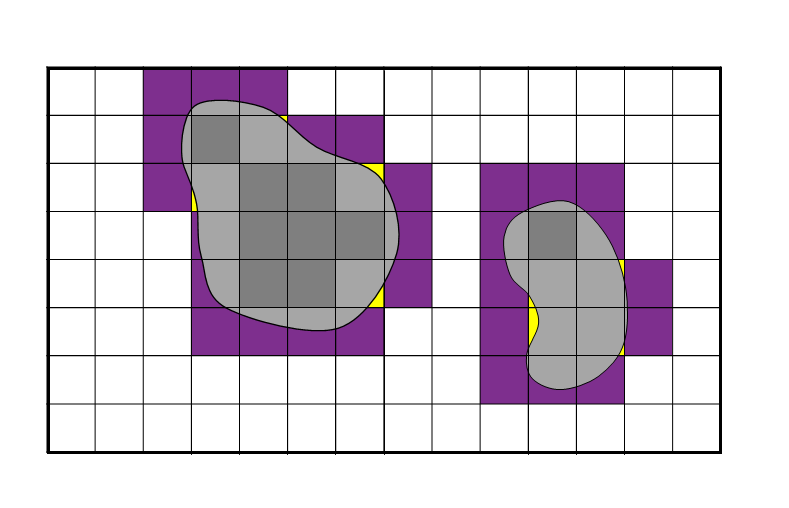}
    \caption{An example of a Cartesian cut mesh with cut elements shown in purple and yellow. Notice that for an equivalent quad-tri mesh every cut element would be decomposed into \textit{at least} one curvilinear triangular element.}
    \label{fig:cut-mesh}
\end{figure}

As high-order elements have more DOF per element than lower order, to maintain reasonable computational cost they are best used on relatively coarse meshes, where ``relatively" is interpreted with respect to an adequately refined mesh for low-order versions of that method. Fitted triangular meshes can require a large amount of elements to resolve a boundary well, potentially leading to bloated computational costs regardless of order. Here, we instead use Cartesian cut meshes, which feature a simple, unfitted Cartesian background mesh ``cut" by the embedded boundaries to yield a hybrid mesh of Cartesian and cut elements as shown in Figure \ref{fig:cut-mesh}.

Cut meshes have a long history in fluid dynamics problems, with early work in the 1970's \cite{reed-originalCutCell}. Cut meshes allow embedded boundaries to be represented with high accuracy while simultaneously maintaining a coarse mesh away from embedded boundaries. The cost of flux-differencing methods is driven by the number of evaluations of the numerical fluxes, which scales linearly with the number of faces in the mesh and nonlinearly with the polynomial degree of the solution. Cut meshes feature fewer faces and thus fewer evaluations of the numerical flux functions than an equivalent quad-tri mesh and (typically) fewer elements and faces than fitted triangular meshes, making them a computationally economic partner for high-order methods. 

While cut meshes are attractive for use with high-order methods, they present a number of challenges. Namely, cut elements can be arbitrarily shaped and sized leading to potentially extremely restrictive CFL conditions (i.e., the small cell problem) and difficulty in constructing high-quality quadrature rules. High-quality quadrature rules are necessary for SBP operators, which classically require diagonal norm mass matrices among other conditions. 

On non-hybrid meshes (i.e. those consisting of a single element type) the cost of constructing SBP-appropriate quadrature (which can be expensive to compute) and the corresponding SBP operators is mitigated by the ability to construct the operators once on a single reference element for all elements in the mesh. The arbitrary geometries of cut elements, besides requiring custom quadrature rules for every cut element, cannot be mapped to a single reference element and thus greatly increase the difficulty and computational cost of constructing traditional SBP operators on cut meshes. 

Even when limited to known element geometries, SBP operators present a new issue when used on hybrid meshes: the conditions needed to induce the classical diagonal-norm SBP property can differ between different types of elements. For example, despite SBP operators being well-established on simplices and tensor-product elements individually, the differing conditions on each type of element create compatibility issues on hybrid quad-tri meshes \cite{chan-SkewSymmES}. However, relaxations of the SBP property exist, such as the generalized SBP property and hybridized SBP property \cite{delReyFernandez-SBPframework, chan-hybridizedSBP, chan-SkewSymmES}, which reduce the cost of constructing SBP operators.

Of particular interest to us are skew-hybridized SBP operators, which were developed by Chan in \cite{chan-SkewSymmES} to overcome the compatibility issue on hybrid quad-tri meshes. Hybridized SBP operators were originally introduced as ``decoupled SBP" operators in \cite{chan-hybridizedSBP} but subsequently became known as ``hybridized" SBP operators \cite{chen-ES-review}. ``Skew" hybridized SBP operators are so named due to their use of the skew-symmetric matrix $(\Q_d - \Q_d^T)$. Given non-SBP differentiation operators $\Q_d$ and boundary integration operators $\B_d$, skew-hybridized SBP operators $\Q_{H,d}$ and $\B_{H,d}$ can be constructed via
\begin{equation}
    \Q_{H,d} = \frac{1}{2}\begin{bmatrix}
        \Q_d - \Q_d^T & \bm{E}^T\B_d \\ 
        -\B_d^T \bm{E} & \B_d
    \end{bmatrix}\label{eq:Q-HSBP}
\end{equation}
and 
\begin{equation}
    \B_{H,d} = \begin{bmatrix}
        \bm{0} & \\ & \B_d
    \end{bmatrix}
\end{equation}
for $d = 1,2$. The skew-hybridized SBP operators satisfy the hybridized SBP property:
\begin{equation}
    \Q_{H,d} + \Q_{H,d}^T = \B_{H,d},
\end{equation}
which is a block version of the traditional SBP property.

The matrix $\bm{E}$ used in \eqref{eq:Q-HSBP} is an extrapolation operator taking the solution values at volume quadrature points and mapping them to values at surface quadrature points (its exact definition will be given later in Section \ref{sec:operators}).  An important property of $\bm{E}$ is that it is dependent on both the volume and surface quadrature rules \cite{chan-hybridizedSBP, chan-SkewSymmES}. With this dependence on volume and surface quadrature, skew-hybridized SBP operators can be thought of as a volume-based derivative operator ($\Q_d$) with boundary-based correction terms coming from the terms with $\bm{E}$ \cite{chan-SkewSymmES, chan-collocation}. 

Importantly, the accuracy of a skew-hybridized SBP operator as a differentiation operator is now tied to the accuracy of the both the surface and volume quadrature \cite{chan-SkewSymmES, chan-hybridizedSBP}. However, since skew-HSBP operators do not require $\Q_d$ to be a SBP operator, they can be readily extended to cut meshes provided sufficiently exact quadrature rules can be constructed on cut elements.

\subsection{Quadrature on Cut Elements}
A number of methods exist for constructing quadrature rules on cut meshes/arbitrary geometries. In some cases, the shape of cut elements is restricted (e.g. to simplicies), allowing established quadrature rules to be applied \cite{stavrev-knownShapeCutCells}. Approaches for constructing quadrature rules typically fall into two categories: generative methods, which numerically create a custom quadrature rule for a specific cut element, and subdivision approaches, which divide the cut element into a collection of subelements who have established quadrature rules (e.g., subtriangulation) from which a composite rule is constructed. 

A common starting point for many generative approaches is moment-fitting, where quadrature rules are constructed on a set of points from target integrals \cite{davis-momentFitting, garhoum-momentFitting1, bui-momentFitting2}. Not all moment-fitting approaches guarantee polynomial exactness (which is often dependent on boundary representation) and positive quadrature weights, though some, such as \cite{garhoum-momentFitting1, legrain-momentFitting}, do but may rely on using a possibly-extreme number of points to do so \cite{piazzon-opt-pruning}. Such methods are also subject to the quality of the sampling points used as they affect the conditioning of the system.

There are however, other generative methods for creating quadrature rules. Work by Saye in \cite{saye-Quadrature, saye-QuadraturePolynomials} provides an algorithm for constructing polynomially exact, positive weight quadrature rules for geometries defined by polynomial level-sets and is commonly used in the cut mesh community. However, Saye's algorithm depends on geometries defined by polynomial level-sets, though the algorithm has also been extended to non-polynoimally defined level-sets via localized projection \cite{kaur-DGDiff}. Here, however, we use explicit parameterizations to represent cut boundaries as we did previously in \cite{taylor-L2-cutDG} and thus are not eligible to use Saye's routines.

Subdivision approaches use existing quadrature on known geometries, such as simplicies and tensor-product elements, to construct a composite quadrature rule on a single cut element. Examples include quad/oct-tree quadrature \cite{kudela-quadtreeQuadr, stavrev-knownShapeCutCells} and subtriangulation approaches \cite{burman-cutFEM1, burman-cutFEM2}. While subdivision approaches have the advantage of leveraging known quadrature formulas and efficient, well-developed libraries, they--by construction--contain many more quadrature points than is needed for a given polynomial degree by virtue of the subdivision. In the case of moving boundaries, the modularity of subdivision-based rules can be beneficial, however, for static boundaries the cost is less permissible.

While the skew-hybridized SBP formulation only requires polynomially exact, positive weight quadrature, the number of quadrature points is an important factor in the computational costs of a scheme both in runtime and memory. Therefore rules using an excessive number of points, such as non-negative least-squares and subdivision, are less preferable with respect to computational cost. However, schemes yielding/using a large number of points can also be quite robust and utilize well-establish libraries and thus warrant further consideration. 

\subsubsection{Quadrature Pruning Approaches}
In addition to routines which provide quadrature rules, there also exist quadrature \textit{pruning} algorithms, which take an initial, many point quadrature rule and prune it to create a quadrature rule on a relatively small/efficient subset of the original points. Importantly, there are pruning algorithms which maintain positivity/non-negativity of the quadrature weights and polynomial exactness such as in \cite{hernandez-quadrature, slobodkins-quadrature-nodeElimination}. Pruning routines can be optimization-based such as in \cite{piazzon-opt-pruning} or iterative processes \cite{vandenbos-caratheodory, wilson-caratheodory}. However, like non-negative moment fitting, pruning routines are generally expensive computationally. In the case of static meshes, however, this cost is only paid once as a preprocessing step.

Pairing a pruning algorithm with an excessive-point quadrature rule allows robust and versatile subdivision-based (and other robust many-point quadrature routines) to be used without excessive time and memory costs during time integration. Exactly how few points a quadrature rule can be reduced to depends on the geometry and the polynomial degree (when polynomial exactness is required). 

Carath\'eodory's theorem, originally developed in the context of convex geometry \cite{caratheodory-thm}, can provide a useful upper bound on the minimum number of quadrature points needed. Theorem \ref{thm:caratheodory-orig} states Cara\'eodory's theorem in its original form while Theorem \ref{thm:caratheodory-quad} states it as applied to quadrature rules. 
\begin{theorem}
    \textbf{(\Caratheodory's Theorem \cite{caratheodory-thm})} Let $S \subset \realR^m$. Any point $\bm{p} \in \conv(S)$ (the convex hull of $S$), can be expressed as a convex combination of at most $(m+1)$ points in $S$.\label{thm:caratheodory-orig}
\end{theorem}

\begin{theorem}
    \textbf{(\Caratheodory's Theorem as applied to quadrature)} Let $D \subset \realR^\Nspace$ be a set of non-zero measure and $N^* = \dim(\totalDegree^N(D))$. Given any $M$-point, non-negative weight quadrature rule exact for $\forall p \in \totalDegree^N(D)$, we can generate a $(N^*+1)$-point non-negative weight quadrature rule exact $\forall p \in \totalDegree^N(D)$ whose quadrature points are a subset of the original rule's points.\label{thm:caratheodory-quad}
\end{theorem}

The proof of Theorem \ref{thm:caratheodory-quad} given Theorem \ref{thm:caratheodory-orig} is very straightforward: the points $\bm{p} \in \realR^m$ are replaced by vectors of each polynomial basis function evaluated at every quadrature point and the coefficients of the convex combination correspond to the new, normalized quadrature weights. The final quadrature weights must remove the scaling from the normalization. While 

Carath\'eodory's theorem provides a theoretically-motivated (when exact quadrature is desired) stopping point for quadrature pruning techniques, but it does not provide any means for obtaining the reduced/pruned quadrature rule. A number of quadrature pruning techniques exist based on various techniques, but a common factor is the use of null vectors of the transpose of the quadrature points' Vandermonde matrix. Such null vectors can be used to manipulate the quadrature weights without altering polynomial exactness.

Here, to generate polynomially exact, positive/non-negative weight quadrature on cut elements we pair subtriangulation with a pruning technique introduced by van den Bos in \cite{vandenbos-caratheodory}, which computes the necessary null vectors via QR factorization. We terminate the pruning according to the upper bound on quadrature points given by Carath\'eodory's theorem. The combination of a pruning routine with the upper bound from Carath\'eodory's theorem we call Carath\'eodory pruning. Such an approach yields efficient, exact, non-negative quadrature on cut elements of arbitrary shape. There are other, potentially more efficient and theoretically attractive (e.g., better controlling the distribution of the quadrature weights) means of generating exact, strictly positive quadrature, such as Vioreanu quadrature \cite{vioreanu-quadrature}, however, for our purposes the combination of subdivision and Carath\'eodory pruning is sufficient.

The combination of non-negative-weight, exact quadrature, the skew-hybridized SBP formulation of \cite{chan-SkewSymmES}, and entropy conservative/stable fluxes constitute, to the authors' knowledge, the first instance of a high-order accurate entropy stable method on cut meshes; this is the major contribution of this work to the literature. The code for our work is written in the programming language Julia and the routines developed to generate our cut meshes and operators are provided in the Julia packages \texttt{PathIntersections.jl} \cite{Julia-PathIntersections} and \texttt{StartUpDG.jl} \cite{Julia-StartUpDG}, respectively. Additionally the codes used for the experiments in this paper are provided in the repository \cite{code-repo}.

While addressing the small cell problem is not one of the objectives of this paper, we do note that state redistribution, first introduced by Berger and Giuliani in \cite{berger-stateRedistr} and extended to high-order methods in \cite{giuliani-DG-SRD}, is an extremely robust \cite{berger-stateRedistr, taylor-L2-cutDG} and powerful method for addressing the small cell problem. While we do not theoretically adapt state redistribution to the entropy-stable setting here, we do present promising numerical results demonstrating it can be used with an entropy stable scheme without numerically violating entropy stability. We use state redistribution as originally presented in \cite{berger-stateRedistr, giuliani-DG-SRD}, however, a newer, more sophisticated and less-dissipative form of state redistribution has since been presented in \cite{giuliani-weightedSRD} and \cite{berger-weightedSRD}. 

In the rest of this paper we will discuss the details of the skew-hybridized SBP formulation, our cut mesh generator, and our Carath\'eodory pruning algorithm in Sections \ref{sec:methods}. In Section \ref{sec:numericalExperiments} we show results for our quadrature generation routine and numerically verify the accuracy and entropy conservation/stability of our scheme. Lastly, we conclude with Section \ref{sec:conclusions} where we discuss future work and open problems.

\section{Methods}\label{sec:methods} 
\subsection{The Skew-Hybridized Summation-by-Parts Formulation}\label{sec:skew-HSBP}
\subsubsection{Solution and Test Spaces}
First, we define the solution and test spaces for our formulation. We assume without loss of generality a Cartesian embedding domain $\CartDomain = [-1,1]^2$ over which we impose a uniform $n_x \times n_y$ Cartesian background grid. The boundaries of the embedded regions $\embeddedRegions$ are allowed to cut the background grid to produce a Cartesian cut mesh $\mesh_h$ of $N_h$ non-overlapping elements.

On Cartesian elements, we use a reference element $\refElement = [-1,1]^2$ to define the solution. On the reference element we define the solution and test spaces as $\tensorProduct^N(\refElement)$, the space of tensor product of polynomials of degree $N$. The space $\tensorProduct^N(\refElement)$ is defined by basis functions
\begin{equation}
    \basisp_{ij}(\bx) = \basisp_{ij}(x,y) = \ell_i(x)\ell_j(y), \quad i,j = 0,1,..., N
\end{equation}
where $\ell_i, \ell_j$ are 1D Lagrange polynomials in the $x$ and $y$ directions, respectively.

On cut elements we define the solution and test spaces on the physical element. For a given cut element $\element^k$, we take the solution and test spaces as 
\begin{equation}
    \totalDegree^N(\element^k) = \text{span} \{~ x^i y^j ~:~ i + j \leq N, ~i,j =0,1,..., N ~\} ,
\end{equation}
the space of total degree $N$ polynomials on $\element^k$. We use the 2D Lagrange polynomials basis functions
\begin{equation}
    \basisp_i(\bx) = \basisp_i(x,y) = \ell_i(x,y), \quad \ell_i(x_j, y_j) = \delta_{ij}, \quad i,j=0,1,..., N
\end{equation}
whose nodal points $\{\bx_i\}_{i=0}^N = \{(x_i,y_i)\}_{i=0}^N \subset \element^k$ are defined by approximate Fekete points. Fekete points are a class of interpolation nodes that yield well-conditioned Vandermonde matrices \cite{sommariva-feketeQR}. We generate approximate Fekete points using the same techniques presented in our previous work on an energy stable cut DG method in \cite{taylor-L2-cutDG}.

For ease of notation, on a given element $\element^k$ we denote the solution and test spaces as $\solnSpace^N(\element^k)$, where 
\begin{equation}
    \solnSpace^N(\element^k) = \left\{ \begin{matrix}
        \tensorProduct^N(\refElement), & \element^k ~\text{Cartesian} \\[3pt]
        ~\totalDegree^N(\element^k), & \element^k ~\text{cut}
    \end{matrix} \right. .
\end{equation}

\subsubsection{Operators}\label{sec:operators}
Now we define the operators for the skew-hybridized formulation. Given volume quadrature  $\{(w_{q,i}, \bx_{q,i})\}_{i=1}^{n_q}$ and face quadrature $\{(w_{f,i}, \bx_{f,i})\}_{i=1}^{n_f}$ rules for a given element $\element^k$, we define the interpolation matrices $\bV_q, \bV_f$ as
\begin{align}
    (\bV_q)_{ij} &= \basisp_j\LRp{\bx_{q,i}}, \quad i=1,..., n_q,\\[6pt]
    (\bV_f)_{ij} &= \basisp_j\LRp{\bx_{f,i}}, \quad i=1,..., n_f,
\end{align}
for $j = 1,..., \dim{\solnSpace^N(\element^k})$ where $\bx_{q,i}$ is $i^{th}$ volume quadrature point, $\bx_{f,i}$ the $i^{th}$ face quadrature point, and $n_q, n_f$ the number of volume and face quadrature points, respectively. Define the matrix of volume quadrature weights $\bW = \diag\{w_{q,1},..., w_{q,n_q}\}$ and the similarly the matrix $\bW_f = \diag\{w_{f,1},..., w_{f,n_f}\}$ of the face quadrature weights. The mass matrix $\M$ is then given by 
\begin{equation}
    \M = \bV_q^T \bW \bV_q
\end{equation}
and the boundary integration matrix in the $d^{th}$ direction $\B_d$ as
\begin{equation}
    \B_d = \bW_f \diag \{\bn_f^{(d)}\}
\end{equation}
where $\bn_f^{(d)}$ if the vector of the $d^{th}$ component of the outward normal at evaluated at every face quadrature point. Note that the boundary integration matrix is diagonal and thus $\B_d = \B_d^T$. We assume that the volume quadrature rule is sufficiently accurate for the mass-matrix to be positive definite and thus invertible (we will explore the conditions each quadrature rule must satisfy in greater detail in Section \ref{sec:quadrature-conditions}).

We define the differentiation (with respect to $\bx_d$) matrix, $\bm{D}_d$, via the relation 
\begin{equation}
    \pderiv{\basisp_i}{\bx_d} = \sum_{j=1}^{\dim (\solnSpace^N(\element^k)} (\D_d)_ij \basisp_j(\bx).
\end{equation}

Next, we define the quadrature-based projection matrix $\bP_q$ as
\begin{equation}
    \bP_q = \M^{-1}\bV_q^T \bW.
\end{equation}
The quadrature-based projection operator projects an arbitrary function evaluated at the volume quadrature points of $\element^k$ onto $\solnSpace^N(\element^k)$ and is vital for the entropy conservation/stability of the scheme. For conserved variables $\bu_h \in \solnSpace^N(\element^k)$, $\bP_q$ is used to project the entropy variable $\bv(\bu_h)$ onto $\solnSpace^N(\element^k)$ via the $L_2$-based projection
\begin{equation}
    \int_{\element^k} (\Pi p)q = \int_{\element^k} pq, \quad \forall q \in \solnSpace^N(\element^k).
\end{equation}
Using $\bP_q$, we can define the \textit{entropy-projected} conserved variables, $\uAux$, as
\begin{equation}
    \uAux = \bu \LRp{\begin{bmatrix}
        \bV_q \\ \bV_f 
    \end{bmatrix}\bP_q \bv\LRp{\bV_q \bu_h} },
\end{equation}
which are defined at all volume and face quadrature nodes. We denote $\uAux$ evaluated at just the face quadrature points as $\uAux_f$.

Given $\bP_q$, we can also define the extrapolation matrix
\begin{equation}
    \E = \bV_f \bP_q
\end{equation}
which extrapolates a function evaluated at the volume quadrature nodes to the face quadrature nodes via its polynomial projection from $\bP_q$. 

With all these operators defined, we can define $\Q_d$, the \textit{generalized} SBP \cite{chan-SkewSymmES} operator in the $d^{th}$ direction, as
\begin{equation}
    \Q_d = \bP_q^T \M \D_d\bP_q.
\end{equation}
Notice the generalized SBP operator and the boundary integration operator do not satisfy the typical SBP property, which would require
\begin{equation}
    \Q_d + \Q_d^T = \B_d,
\end{equation}
which is not in fact even dimensionally consistent for $\Q_d$ and $\B_d$ as defined above. However, from the generalized SBP operator, we can define the skew-hybridized SBP operator $\Q_{H,d}$:
\begin{equation}
    \Q_{H,d} = \frac{1}{2} \begin{bmatrix}
        \Q_d - \Q_d^T & \E^T\B_d \\ -\B_d \E & \B_d
    \end{bmatrix}.
\end{equation}
The skew-hybridized SBP operator then satisfies the hybridized SBP property:
\begin{equation}
    \Q_{H,d} + \Q_{H,d}^T = \begin{bmatrix}
        \bm{0} & \\ & \B_d
    \end{bmatrix}
\end{equation}
for sufficiently accurate surface and volume quadrature \cite{chan-SkewSymmES}. The exact conditions on the quadrature rules will be discussed in Section \ref{sec:quadrature-conditions}.

\subsubsection{Formulation}
Given the operators defined in Section \ref{sec:operators} and entropy conservative two-point fluxes $\fEC^d, d=1,2$, the skew-hybridized formulation on a single element $\element^k$ for semi-discrete solution $\bu_h \in \solnSpace^N(\element^k)$ is given by:
\begin{equation}
    \M \deriv{\bu_h}{t} + \sum_{d=1}^{2} \begin{bmatrix}
        \bV_q \\ \bV_f
    \end{bmatrix}^T 2\LRp{\Q_{H,d} \circ \F_d}\ones + \bV_f^T \B_d \LRp{\f^*_d - \f_d(\uAux_f} = \bm{0}\label{eq:sHSBP-formulation}
\end{equation}
where $\ones$ is an appropriately sized vector of all ones and the flux matrices $\F_d$ are given by
\begin{equation}
    (\F_d)_{ij} = \fEC^d(\uAux_i, \uAux_j), \quad i,j = 1,..., (n_q + n_f).\label{eq:fluxMatrices}
\end{equation}

The vectors $\f^*_d, d= 1,2$ denote the boundary fluxes and can be defined using either the same two-point, entropy-conservative fluxes used in the volume, $\fEC^d, d=1,2$, or entropy-stable fluxes $\fES^d, d=1,2$. For an entropy conservative scheme $\f^*_d$ is given by
\begin{equation}
    (\f^*_d)_i := \fEC^d\LRp{(\uAux_f^+)_i, (\uAux_f)_i}, \quad i = 1,..., n_f
\end{equation}
where $\uAux_f^+$ is $\uAux_f$ on the appropriate face-neighbor of $\element^k$. For an entropy stable scheme $\f^*_d$ is instead defined by the three-argument entropy stable flux
\begin{equation}
    \bn^{(d)}_i (\f^*_d)_i := \fES^d\LRp{(\uAux_f^+)_i, (\uAux_f)_i, \bn^{(d)}_i}
\end{equation}
for $i=1,..., n_f$ where $\bn^{(d)}_i$ is the $d^{th}$ component of the outward normal vector at the $i^{th}$ face quadrature point. Note for entropy stable fluxes defined as shown in Equation \eqref{eq:ES-ds-flux} this definition is necessary to avoid scaling the dissipation term in the entropy stable flux by the outward normal:
\begin{equation}
    \bn^{(d)}_i (\f^*_d)_i = \bn^{(d)}_i \fEC^d\LRp{(\uAux_f^+)_i, (\uAux_f)_i} + \bm{d}_{\entropy}^d\LRp{(\uAux_f^+)_i, (\uAux_f)_i}.
\end{equation}

As shown in \cite{chan-SkewSymmES, chan-hybridizedSBP}, this formulation is high-order accurate, entropy conservative for sufficiently accurate quadrature, and locally conservative (which is necessary to show convergence under mesh refinement to the weak solution via a Lax-Wendroff type argument as defined in \cite{shi-shu-localConservation}). Next we discuss the exact conditions on quadrature needed for \eqref{eq:sHSBP-formulation} to be high-order accurate and entropy conservative/stable.

\subsubsection{Conditions on Quadrature/The Skew-Hybridized SBP Operator}\label{sec:quadrature-conditions}
For convenience, we again note the definition of the skew-hybridized operator:
\begin{equation}
    \Q_{H,d} = \frac{1}{2} \begin{bmatrix}
        \Q_d - \Q_d^T & \E^T\B_d \\ -\B_d \E & \B_d.
    \end{bmatrix}
\end{equation}
which by construction satisfies the hybridized SBP property:
\begin{equation}
    \Q_{H,d} + \Q_{H,d}^T = \begin{bmatrix}
        \bm{0} & \\ & \B_d
    \end{bmatrix}.
\end{equation} 
$\Q_{H,d}$ can be thought of as a volume-based differentiation operator (from the contribution of $\Q_d$) with boundary-based correction terms from the off-diagonal blocks involving $\E$ and $\B_d$ \cite{chan-SkewSymmES}. A consequence of its dependence on $\E$  and $\B_d$, however, is that the accuracy of $\Q_{H,d}$ as a differentiation operator is now tied to the accuracy of the volume (via $\E$) and surface (via $\B_d$) quadrature rules.

As shown in \cite{chan-hybridizedSBP}, the formulation given in \eqref{eq:sHSBP-formulation} is entropy stable/conservative if the skew-hybridized SBP operator $\Q_{H,d}$ can satisfy
\begin{equation}
    \Q_{H,d} \ones = \bm{0}.\label{eq:diff-constants}
\end{equation}
Given $\Q_{H,d}$ acts a differentiation matrix, Condition \eqref{eq:diff-constants} corresponds to the ability to exactly differentiate constant functions \cite{chan-SkewSymmES}.

While \eqref{eq:diff-constants} is needed to ensure entropy stability/conservation, it is not sufficient to ensure high-order accuracy. To maintain high-order accuracy $\Q_{H,d}$ must be able to exactly differentiate up to degree $2N-1$ polynomials (this comes from the need to exactly integrate the product of a degree $N$ polynomial and the derivative of a degree $N$ polynomial). As shown in \cite{chan-SkewSymmES}, this condition for high-order accuracy requires the volume quadrature rule used to construct $\E$ to be exact for polynomials of degree $2N-1$ and the surface quadrature rule used for $\B_d$ to be exact for polynomials of degree $2N$. 

In addition to the conditions imposed by $\Q_{H,d}$ on the quadrature rules, entropy stability/conservation also requires the mass-matrix to be positive-definite. While the mass matrix is \textit{guaranteed} to be positive definite (and exact) if the volume quadrature rule is exact for polynomials of degree $2N$, the mass matrix will remain positive-definite even for inexact volume quadrature if the quadrature weights are non-negative and a sufficient number of the weights are strictly positive; these conditions are easily satisfied by the same conditions needed for $\Q_{H,d}$. 

On Cartesian elements, the above conditions on quadrature are easily satisfied. On cut elements, we appeal to subtriangulation and Carath\'eodory pruning to produce the necessary quadrature rules, which we discuss next.

\subsection{Cut Mesh Generation and Quadrature on Cut Elements}\label{sec:caratheodryPruning}
We construct our cut meshes using the Julia libraries \texttt{PathIntersection.jl} and \texttt{StartUpDG.jl}. Our cut meshes are constructed using explicit parameterizations of the embedded boundaries as previously detailed in \cite{taylor-L2-cutDG}.
Using explicit parameterization has a number of benefits, such as allowing piecewise-defined curves and thus an easy interface for imposing boundary conditions on portions of the cut boundaries and exact representation of sharp features. In particular, our code uses \textit{stop points} to define both user-defined subintervals and junctions in piecewise curves. 

\begin{figure}[H]
    \centering
    \includegraphics[width=0.75\linewidth]{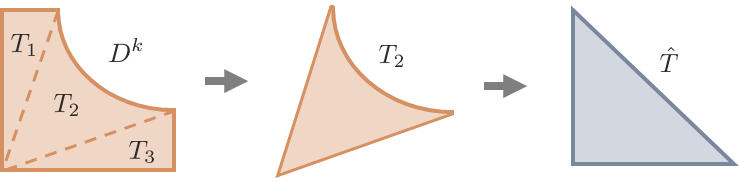}
    \caption{An example of the subtriangulation process used on each cut element mapping a subtraingle to the reference triangle $\hat{T}$.}
    \label{fig:subtriangulation}
\end{figure}

Our code provides a piecewise representation of the boundary of every cut element, with stop points denoting the endpoints of the element's faces. We use the stop points to triangulate each cut element via the Julia package \texttt{Triangulate.jl}. For a given cut element $\element^k$, we decompose it into $N_T$ polynomially-defined non-overlapping curvilinear triangles $T_i$:
\begin{equation}
    \element^k = \bigcup_{i=1}^{N_T} T_i
\end{equation}
For each triangle, we construct a geometric mapping $\bm{g}_d:\hat{T} \to T_i, ~d= 1,2$ from a reference triangle, $\hat{T}$, to each subtriangle $T_i$. Note for high-order accuracy the polynomial degree of the geometric mappings must be isogeometric, i.e., at least the same degree as the degree of the solution. We take the mappings to be total degree $N$ in all meshes, i.e. $\bm{g}_d \in \totalDegree^N(\hat{T})$. An example of the subtriangulation process is shown in Figure \ref{fig:subtriangulation}.

The geometric mappings provide both a polynomial approximation of the typically non-linear curved boundaries and access to a reference element, allowing both exact surface and volume quadrature rules to be constructed. However, the geometric mapping and its Jacobian also increase the degree of the integrand. 

We will illustrate the effect of the geometric mapping with respect to the volume integrals. Consider a polynomial $p \in \totalDegree^M(\element^k)$. On a given subtriangle $T_i \subseteq \element^k$, to map an integral with respect to $\bx \in T_i$ to an integral with respect to $\bhat{x} \in \hat{T}$, we take $\bx_d = \bm{g}_d
(\bhat{x})$, $d=1,2$. The subintegral on $T_i$ is then related to its corresponding integral on $\hat{T}$ by 
\begin{equation}
    \int_{T_i} p(\bx) \dint \bx = \int_{\hat{T}} p(\bm{g}(\bhat{x})) \LRb{\pderiv{\bm{g}}{\bhat{x}} }\dint \bhat{x} = \int_{\hat{T}} \hat{p}(\bhat{x}) \dint \bhat{x}.
\end{equation}
For $p \in \totalDegree^M(\element^k)$ and $\bm{g}_d \in \totalDegree^N(\hat{T}), ~d=1,2$ we have that $p \circ \bm{g} \in \totalDegree^{MN}(\hat{T})$. Similarly, the determinant of the Jacobian $\LRb{\nicefrac{\partial\bm{g}}{\partial\bhat{x}} } \in \totalDegree^{2(N-1)}$. Thus the final integrand on the reference triangle is a polynomial of total degree $MN + 2(N-1)$. For the accuracy conditions on $\Q_{H,d}$, we need $M = 2N-1$, resulting in the final integrand $\hat{p}$ on the reference triangle being total degree $(2N^2 + N -2)$; a quadratic increase in degree with respect to $N$, the polynomial degree of the solution. To exactly integrate the mass matrix would result in the final integrand being degree $2(N^2 + N - 1)$. Table \ref{tab:ref-integrand-deg} shows the polynomial degree of the final integrand on the reference triangle for various values of $N$ for reference.

\begin{table}[H]
    \centering
    \caption{Total degree of the final integrand $\hat{p}$ on $\hat{T}$ after applying the geometric mapping.}
    \label{tab:ref-integrand-deg}
    \begin{tabular}{c|cccccccc}
         $M \backslash N$ & 1 & 2 & 3 & 4 & 5 & 6 & 7 & 8 \\
         \hline\hline
         $2N-1$ &  1 & 8 & 19 & 34 & 53 & 76 & 103 & 134 \\
         $2N$ & 2 & 10 & 22 & 38 & 58 & 82 & 110 & 142 \\ 
         \hline\hline
    \end{tabular}
\end{table}

Table \ref{tab:ref-integrand-deg} helps to illustrate why subdivision type quadrature approaches--when used by themselves--are so inefficient when paired with high-order methods: the subquadrature rules must be constructed for polynomial degrees far in excess of the degree of the solution. Moreover, the cost of these many point subquadrature rules is multiplied by the number of triangles in the triangulation of each cut element.

However, subdivision based quadrature rules are good starting point for a pruning routine: by construction the subquadrature rules are exact and have purely positive weights. We use the QR factorization approach of van den Bos \cite{vandenbos-caratheodory} for pruning, which we describe in Appendix \ref{ap:pruning}.

\subsection{State Redistribution and Entropy Conservation/Stability}
While not the focus of this paper, one of the greatest challenges of/deterrents from using cut meshes is the small cell problem: cut elements can be arbitrarily small/skewed. When using explicit time integration schemes, the small cell problem can result in extremely prohibitive CFL conditions. One method for addressing the small cell problem is state redistribution, which was introduced \cite{berger-stateRedistr} and adapted to high-order methods in \cite{giuliani-DG-SRD}. 

State redistribution is high-order accurate and conserves the average of the solution. State redistribution involves projections spanning multiple elements and averaging of those projections (or more generally convex combinations of those projections \cite{giuliani-weightedSRD, berger-weightedSRD}) and as a result is dissipative in nature. In the context of linear symmetric hyperbolic conservation laws, state redistribution can be added to a $L_2$ energy stable scheme without damaging the scheme's energy stability/conservation as shown in \cite{taylor-L2-cutDG}. However, state redistribution readily violates conservation of entropy when applied to an entropy conservative scheme. 

Encouragingly, in the entropy \textit{stable} case we have observed that state redistribution does not violate the entropy stability of our scheme. We show the effect of state redistribution on a system's entropy in Section \ref{sec:entropy-status}. In keeping with our observations, many of our entropy stable numerical experiment use state redistribution. However, we acknowledge that it remains to be proven whether state redistribution can in fact violate entropy stability; we save this investigation for future work.

\section{Numerical Experiments}\label{sec:numericalExperiments} 

\subsection{Equations, Fluxes, and Time Integration}
In the following numerical experiments we consider two hyperbolic systems, the shallow water equations and the compressible Euler equations, on 2D domains. The definitions of these systems and their respective fluxes and the time integration scheme used are given in this section for reference. 

In both systems, when an entropy stable flux is needed we use the Lax-Friedrichs flux
\begin{equation}
    \fES^d(\bu_L, \bu_R, \bn_d) = \bn_d\avg{\bu} - \frac{\lambda}{2}\jump{\bu}
\end{equation}
where $\lambda$ is the maximum wavespeed or an estimate thereof as previously mentioned. We use the Lax-Friedrichs fluxes provided in the \texttt{Trixi.jl} library, which use the wavespeed approximation by Davis \cite{davis-wavespeed}.

\subsubsection{Shallow Water Equations}
The 2D shallow water equations with constant bathymetry (i.e. bottom height) are given by:
\begin{equation}
    \pderiv{}{t}{}\begin{bmatrix}
        h \\ hu_1 \\ hu_2
    \end{bmatrix} + \pderiv{}{x}{}\begin{bmatrix}
        hu_1 \\ hu_1^2 + \frac{1}{2} gh^2 \\ hu_1 u_2
    \end{bmatrix}  + \pderiv{}{y}{}\begin{bmatrix}
        hu_2 \\ hu_1u_2 \\ hu_2^2 + \frac{1}{2}gh^2
    \end{bmatrix} = \bm{0}\label{eq:SWE}
\end{equation}
for conserved variables $\bu = (h, hu_1, hu_2)^T$ where $h$ is the water height, $u_1$ is the velocity in the $x$-direction, and $u_2$ is the velocity in the $y$-direction. For the entropy conservative flux we use the Wintermeyer \cite{wintermeyer-SWE-flux} flux provided in \texttt{Trixi.jl}.

\subsubsection{Compressible Euler Equations}
The compressible Euler equations are given by:
\begin{equation}
    \pderiv{}{t}{}\begin{bmatrix}
        \rho \\ \rho u_1 \\ \rho u_2 \\ E
    \end{bmatrix} + \pderiv{}{x}{}\begin{bmatrix}
        \rho u_1 \\ \rho u_1^2 + p \\ \rho u_1 u_2\\ u_1(E+p)
    \end{bmatrix}  + \pderiv{}{y}{}\begin{bmatrix}
        \rho u_2 \\ \rho u_1 u_2 \\ \rho u_2^2 + p \\ u_2(E+p)
    \end{bmatrix} = \bm{0}\label{eq:SWE}
\end{equation}
for conserved variables $\bu = (\rho, \rho u_1, \rho u_2, E)^T$ where $\rho$ is density, $u_1$ is the velocity in the $x$-direction, $u_2$ is the velocity in the $y$-direction, $E$ is the total energy, and $p$ is the pressure. The total energy is related to the other variables by:
\begin{equation}
    E = \frac{1}{2}\rho (u_1^2 + u_2^2) + \frac{p}{\gamma - 1}
\end{equation}
where $\gamma$ is the ratio of specific heats for the gas.

We use the entropy conservative flux of Ranocha \cite{ranocha-EC-eulerFluxes, ranocha-EulerFlux} as provided in the \texttt{Trixi.jl} library.

\subsubsection{Time Integration}
In all experiments we use the explicit, adaptive time integration scheme \texttt{Tsit5} \cite{tsitouras-Tsit5} as provided in the Julia package \texttt{OrdinaryDiffEq.jl} \cite{Julia-OrdinaryDiff}. The initial time step used is given in the description of each experiment.

\subsection{Entropy Conservation and Stability}\label{sec:entropy-status}
Here we numerically verify the entropy status of our scheme. For this experiment, we are interested in the evolution of the system's entropy over time. This can be done by considering the integral over the domain of the time derivative of entropy which, for entropy conservative boundary conditions, satisfies
\begin{equation}
    \int_{\domain} \pderiv{\entropy(\bu)}{t} \dint \bx \leq 0
\end{equation}
with equality holding in the absence of shocks. At the semi-discrete level, this entity is called the \textit{entropy residual} and is given by:
\begin{equation}
    \sum_{k=1}^{N_h} \ones^T \bW^{(k)}\pderiv{\entropy(\bu_h^{(k)})}{t}{} = \sum_{k=1}^{N_h} \LRp{\bv_h^{(k)}}^T \M^{(k)} \deriv{\bu_h^{(k)}}{t} \leq 0
\end{equation}
where the superscript $(k)$ indicates the appropriate entity on the $k^{th}$ element and
\begin{equation}
    \bv_h^{(k)} = \bP_q^{(k)} \bv\LRp{\bV_q^{(k)} \bu_h^{(k)} }.\label{eq:entropy-res}
\end{equation}
In the semi-discrete case, equality (up to round-off error) holds if the scheme is entropy conservative, while the inequality hold for an entropy stable scheme.

We first consider the shallow water equations on a Cartesian embedding domain $\CartDomain = [-1,1]^2$ with a circle of radius $R = 0.331$ removed from the center. We impose a $16 \times 16$ background Cartesian mesh. To test the entropy status of our scheme we compute the entropy residual as shown in \eqref{eq:entropy-res} at each time step. We use an initial time step of $\Delta t = 1 \times 10^{-4}$ and integrate up to an end time of $t=1$. 

\begin{figure}[H]
    \centering
    \includegraphics[width=0.75\linewidth]{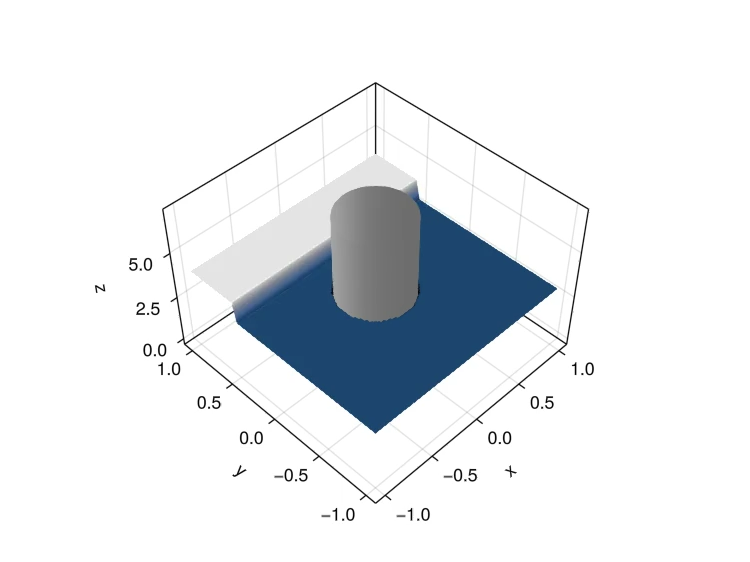}
    \caption{The initial water height for the entropy residual experiments.}
    \label{fig:entropy-res-IC}
\end{figure}

We enforce reflective wall boundary conditions, which are entropy conservative \cite{svard-ES-eulerBCs}, on all boundaries including the circle walls. We start with a zero-velocity solution with discontinuous water height given by
\begin{equation}
    h(x,y,t=0) = \left\{ \begin{matrix}
        3, & \quad y \geq 0.5 \\
        2, & \quad y < 0.5
    \end{matrix}\right. .
\end{equation}
Note that in the continuous setting, weak solutions should dissipate entropy thanks to the discontinuity in the water height. This serves as a stress test for an entropy conservative scheme: it should conserve entropy despite the discontinuity, which in non-entropy conservative/stable scheme can lead to instability. In the entropy stable case, entropy should immediately be dissipated due to the discontinuity and the dissipation wane as the solution approaches its steady state. Figure \ref{fig:entropy-res-IC} shows the initial water height.

\begin{figure}[H]
    \centering
    \begin{subfigure}{0.48\textwidth}
        \includegraphics[width=\linewidth]{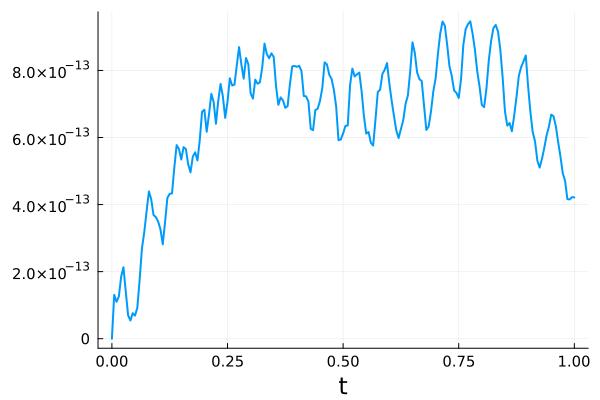}
        \caption{}
        \label{fig:entropy-res-EC-noSRD}
    \end{subfigure}
    \begin{subfigure}{0.48\textwidth}
        \includegraphics[width=\linewidth]{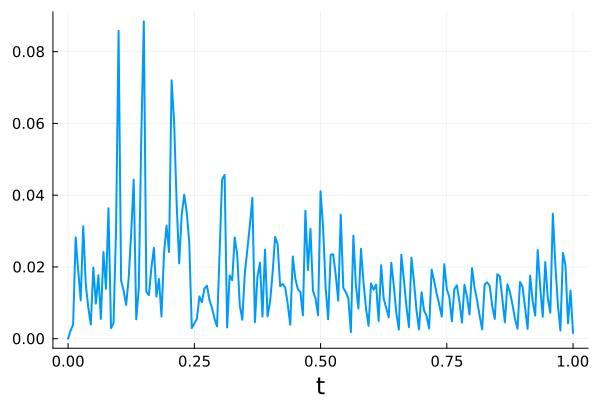}
        \caption{}
        \label{fig:entropy-res-EC-SRD}
    \end{subfigure}
    \caption{The entropy residual for the entropy conservative scheme without (\ref{fig:entropy-res-EC-noSRD}) and with (\ref{fig:entropy-res-EC-SRD}) state redistribution. Without state redistribution, the entropy residual is zero up to numerical round-off error.}
    \label{fig:entropy-res-EC}
\end{figure}

\begin{figure}[H]
    \centering
    \begin{subfigure}{0.48\textwidth}
        \includegraphics[width=\linewidth]{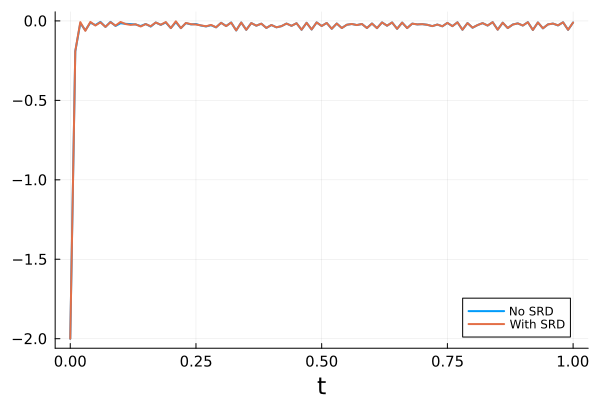}
        \caption{}
        \label{fig:entropy-res-ES}
    \end{subfigure}
    \begin{subfigure}{0.48\textwidth}
        \includegraphics[width=\linewidth]{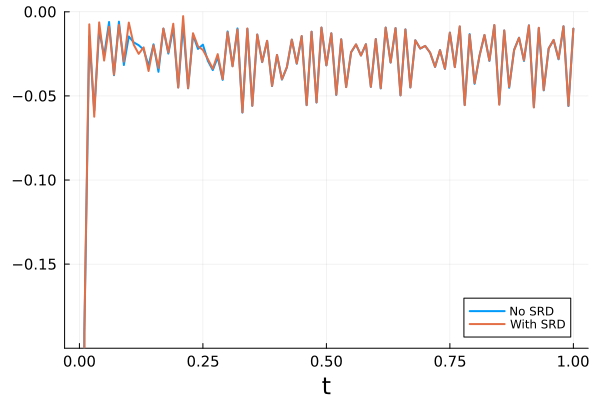}
        \caption{}
        \label{fig:entropy-res-ES_zoomed}
    \end{subfigure}
    \caption{The entropy residual for the entropy stable scheme both with (orange curve) and without (blue curve) state redistribution. As shown in (\ref{fig:entropy-res-ES}), the initial shock produces a large amount of dissipation which then tapers off as the system approaches its steady state. \ref{fig:entropy-res-ES_zoomed} highlights that the entropy residual remains less than or equal to zero both with and without state redistribution. Importantly, the difference between the entropy residual with and without state redistribution is quite small.}
    \label{fig:entropy-res-ES-all}
\end{figure}

We consider four schemes for this experiment: the entropy conservative scheme with and without state redistribution, and the entropy stable scheme with and without state redistribution. All schemes use polynomial degree $N = 4$. Figures \ref{fig:entropy-res-EC} and \ref{fig:entropy-res-ES-all} show the entropy residual for each of the four schemes.

Without state redistribution, the entropy conservative scheme maintains entropy conservation within round-off error despite the presence of the shock. The addition of state redistribution to the entropy conservative scheme however violates entropy conservation by a non-negligible margin. The entropy stable scheme meanwhile maintains entropy stability even when state redistribution is applied. While it is encouraging to see that state redistribution did not violate entropy stability in this instance, we again stress that entropy stability may not always hold under state redistribution.

\subsection{$h$-Convergence Study}
\begin{figure}[H]
    \centering
    \includegraphics[width=0.6\linewidth]{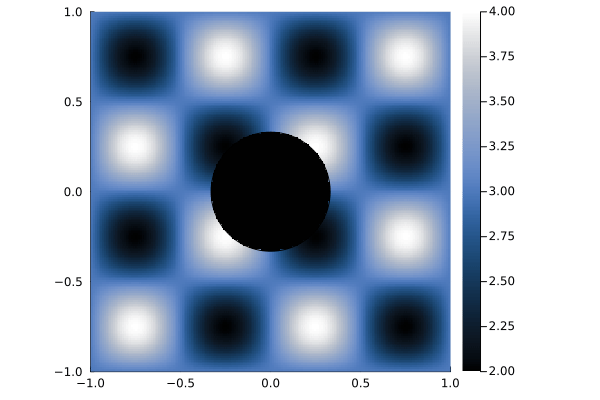}
    \caption{Water height for the manufactured solution at time $t = 0$.}
    \label{fig:h-conv-IC}
\end{figure}

For this experiment we numerically verify the high-order accuracy of our scheme via an $h$-convergence study. We take the Cartesian embedding domain $\CartDomain = [-1,1]^2$ from which we remove a circle of radius $R = 0.331$ from the origin. We consider the manufactured solution to the shallow water equations defined by primitive variables:
\begin{equation}
    \begin{bmatrix}
        h \\ u_1 \\ u_2
    \end{bmatrix} = \begin{bmatrix}
        \sin(2\pi x) \sin(2 \pi y) \cos(\pi t) + 3 \\ 1 \\ 1
    \end{bmatrix}.
\end{equation}
Figure \ref{fig:h-conv-IC} shows the water height at time 0. We start with a $4 \times 4$ Cartesian background mesh and refine in each dimension by a factor of 2 up to a $32 \time 32$ background mesh. We use the manufactured solution values as boundary conditions and simulate up to an end time of $t = 0.3$. We use an initial time step of $\Delta t = 1 \times 10^{-4}$ for the $4 \times 4$ mesh and $\Delta t = 1 \times 10^{-5}$ for the other meshes. 

\begin{figure}[H]
    \centering
    \includegraphics[width=0.5\linewidth]{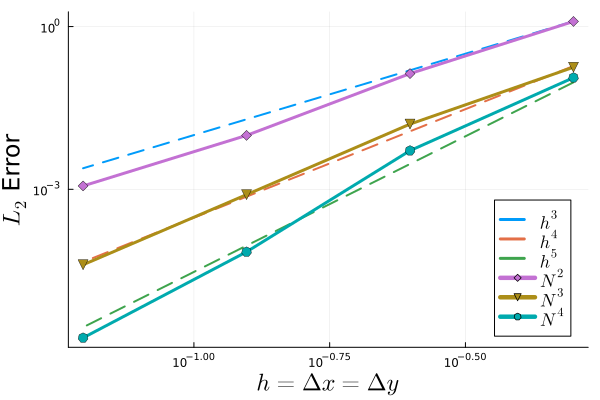}
    \caption{$L_2$ error versus mesh size $h$ for solutions of degree $N = 2, 3, 4$ for the $h$-convergence study.}
    \label{fig:h-conv-plot}
\end{figure}

We conducted the convergence study for solutions of degree $N = 2,3,4$. As shown in Figure \ref{fig:h-conv-plot}, our method recovers the expected convergence rate of $h^{N+1}$ in each case. 

\subsection{Benchmark Problem: Entropy Wave}
Here we consider an analytic solution to the compressible Euler equations. When initialized with constant velocity and pressure, the compressible Euler equations reduce to the advection equation with vector-valued coefficients applied to the initial density field. This special instance of the Euler equations is known as an \textit{entropy wave} and is  one of three modes of linear flow fluctuation in gases (the others being acoustic waves and vorticity waves) \cite{chu-entropyWave}. The reduced equations are:
\begin{equation}
    \pd{\rho}{t}{}\begin{bmatrix} 1 \\ u_1 \\ u_2 \\[6pt] \frac{1}{2}\LRp{u_1^2 + u_2^2} \end{bmatrix} + 
 \pd{\rho}{x}{}\begin{bmatrix} u_1 \\  u_1^2 \\ u_1 u_2 \\[6pt] \frac{1}{2}\LRp{u_1^2 + u_2^2}u_1\end{bmatrix} + 
 \pd{\rho}{y}{}\begin{bmatrix} u_2 \\ u_1 u_2 \\ u_2^2 \\[6pt] \frac{1}{2}\LRp{u_1^2 + u_2^2}u_2\end{bmatrix} = \bm{0}
\end{equation}
which can also be expressed as
\begin{equation}
 \bm{c}_t\pd{\rho}{t}{} + 
 \bm{c}_x\pd{\rho}{x}{} + 
 \bm{c}_y\pd{\rho}{y}{} = \bm{0}. \label{eq:entropy-wave-advection}
\end{equation}

For this experiment, we again consider the embedding domain $\CartDomain = [-1,1]^2$ over which we impose a $16 \times 16$ background Cartesian mesh. A circle of radius $R=0.331$ centered at the origin is removed. We consider the initial condition
\begin{equation}
    p = 3, \quad u_1 = u_2 = \frac{1}{2}, 
\end{equation}
\begin{equation}
    \rho(x,y, t=0) = \sin(2\pi x)\sin(2\pi y) + 2
\end{equation}
which by the method of characteristics applied to the reduced PDE/advection equation \eqref{eq:entropy-wave-advection} yields the analytic solution
\begin{equation}
    \rho(x,y, t) = \sin\LRp{2\pi(x - u_1 t)}\sin\LRp{2\pi(y - u_2 t)} + 2.
\end{equation}

\begin{figure}[H]
    \centering
    \begin{subfigure}[b]{0.4\textwidth}
    \centering
    \includegraphics[width=\linewidth]{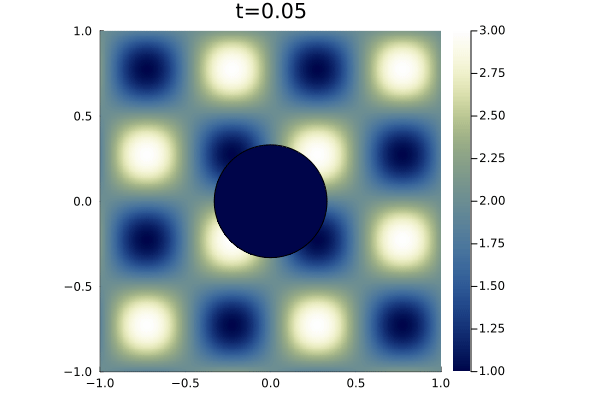}
    \end{subfigure}
    \begin{subfigure}[b]{0.4\textwidth}
    \centering
    \includegraphics[width=\linewidth]{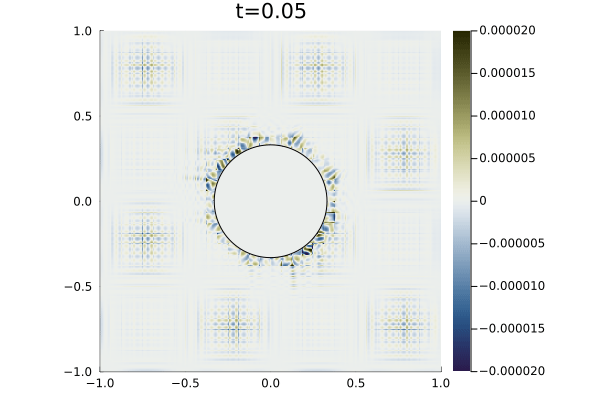}
    \end{subfigure}
    
    \begin{subfigure}[b]{0.4\textwidth}
    \centering
    \includegraphics[width=\linewidth]{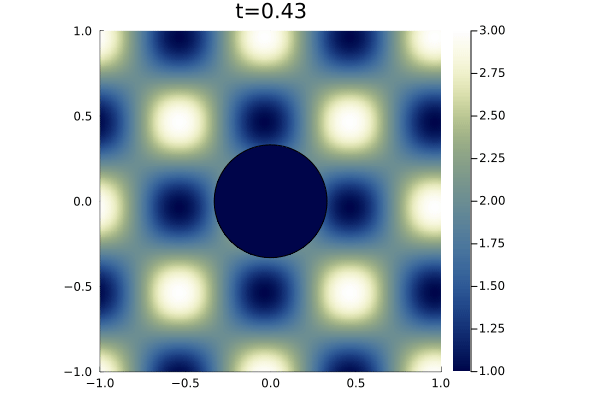}
    \end{subfigure}
    \begin{subfigure}[b]{0.4\textwidth}
    \centering
    \includegraphics[width=\linewidth]{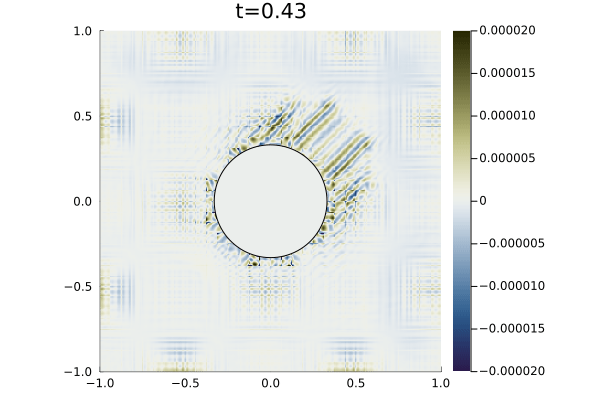}
    \end{subfigure}

    \begin{subfigure}[b]{0.4\textwidth}
    \centering
    \includegraphics[width=\linewidth]{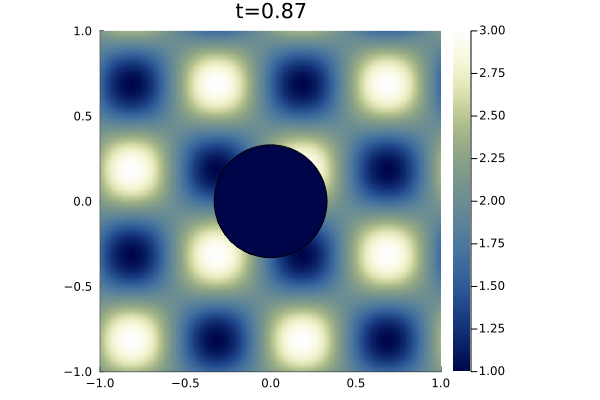}
    \end{subfigure}
    \begin{subfigure}[b]{0.4\textwidth}
    \centering
    \includegraphics[width=\linewidth]{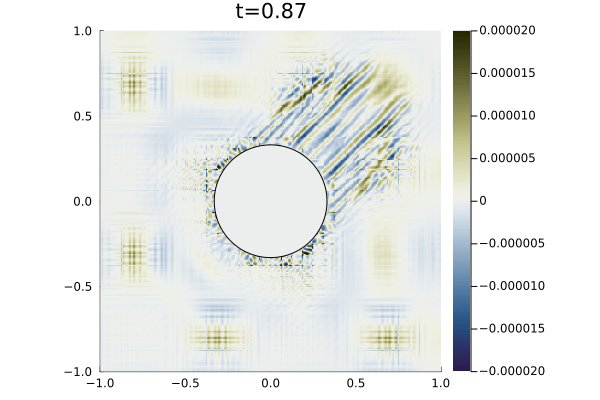}
    \end{subfigure}
    
    \begin{subfigure}[b]{0.4\textwidth}
    \centering
    \includegraphics[width=\linewidth]{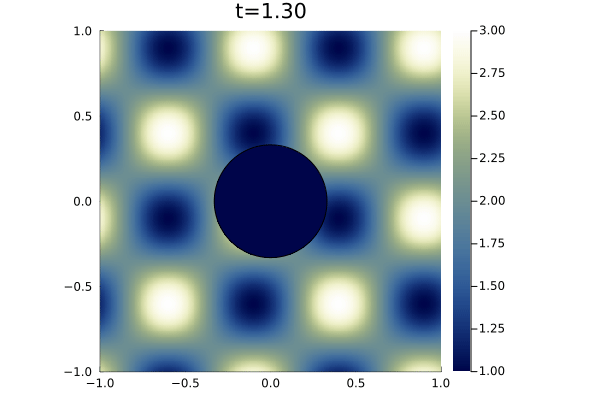}
    \end{subfigure}
    \begin{subfigure}[b]{0.4\textwidth}
    \centering
    \includegraphics[width=\linewidth]{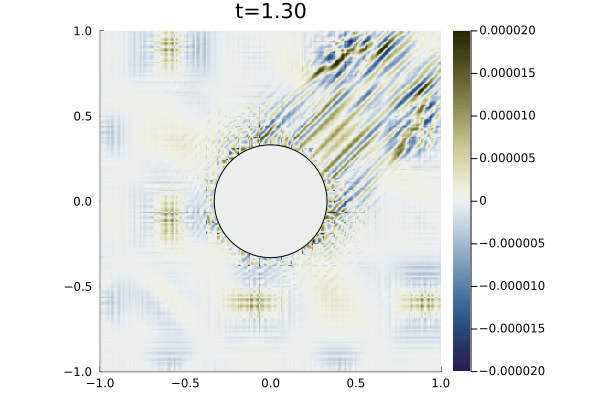}
    \end{subfigure}
    \caption{Snapshots of the density (left column) and error of the density (right column) for the entropy wave simulation.}
    \label{fig:entropyWave-snapshots}
\end{figure}

For boundary conditions, we prescribe the analytic solution on all boundaries. For time integration we use an initial time step of $1\times 10^{-6}$ (a typical timestep was on the order of $1 \times 10^{-4}$) with state redistribution. The simulation was run to an end time of $t=1.3$ with polynomial degree $N = 4$. Figure \ref{fig:entropyWave-snapshots} shows snapshots the solution an error at various points in time alongside the error in the density field. As illustrated in the error plots, cut elements produce higher errors than Cartesian elements, however, the maximum error magnitude is five orders of magnitude less than the magnitude of the solution.

\subsection{Biconvex Supersonic Airfoil}
For this experiment we consider the compressible Euler equations to simulate supersonic flow over a biconvex airfoil. This experiment highlights our method's ability to handle sharp features in the mesh and solution while remaining stable.

\begin{figure}[H]
    \centering
    \includegraphics[width=0.75\linewidth]{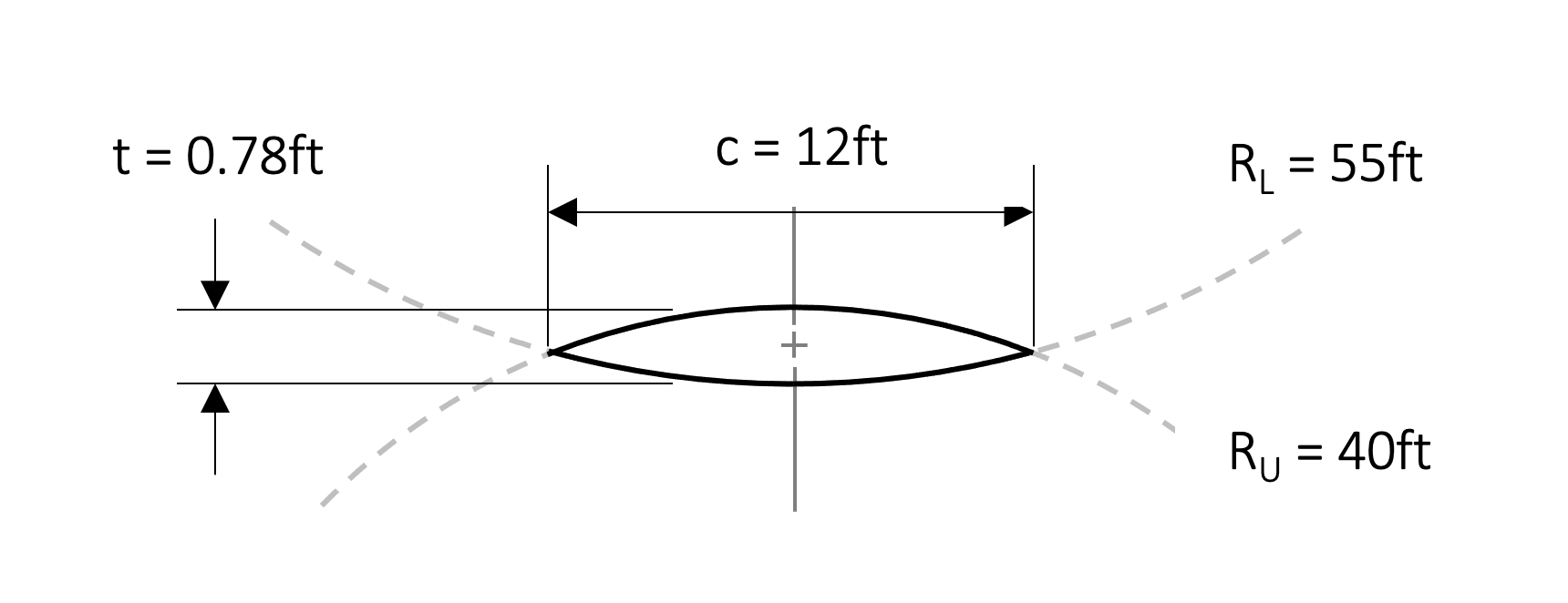}
    \caption{Unscaled geometry of the biconvex airfoil.}
    \label{fig:airfoil-geom}
\end{figure}

For the mesh we take the Cartesian embedding domain $\CartDomain = [-0.5, 5.5] \times [-1.5, 1.5]$ over which we impose a $120 \times 60$ background Cartesian mesh. The airfoil, whose unscaled geometry is given in Figure \ref{fig:airfoil-geom}, is scaled to a non-dimensionalized chord length of 0.5 and centered about the $x$-axis at zero angle of attack and with its chord line at $y = 0$.

We apply reflective wall boundary conditions to the airfoil and the top ($y = 1.5$) and bottom ($y=-1.5$) walls of the domain. We take the left ($x=-0.5$) boundary as the inlet with density $\rho_\infty = 50$, pressure $p_\infty = 50$, $x$-velocity $u_{\infty,1} = $ Mach 1.5, and $y$-velocity $u_{\infty,2} = 0$. Extrapolation boundary conditions were used on the right ($x= 5.5$) boundary.

\begin{figure}[H]
    \centering
    \begin{subfigure}[b]{0.49\textwidth}
        \centering
        \includegraphics[trim={0cm 1cm 1.5cm 0cm}, clip, width=\textwidth]{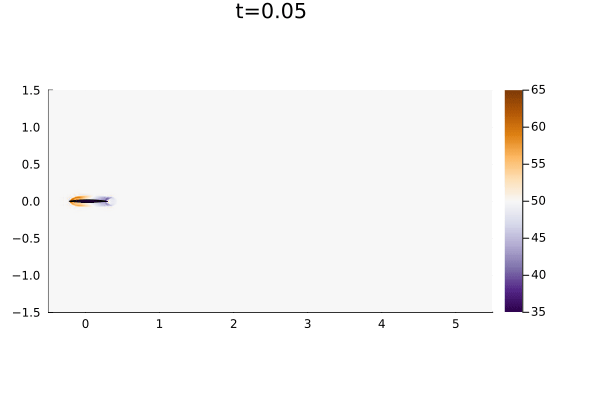}
    \end{subfigure}
    \begin{subfigure}[b]{0.49\textwidth}
        \centering
        \includegraphics[trim={0cm 1cm 1.5cm 0cm}, clip, width=\textwidth]{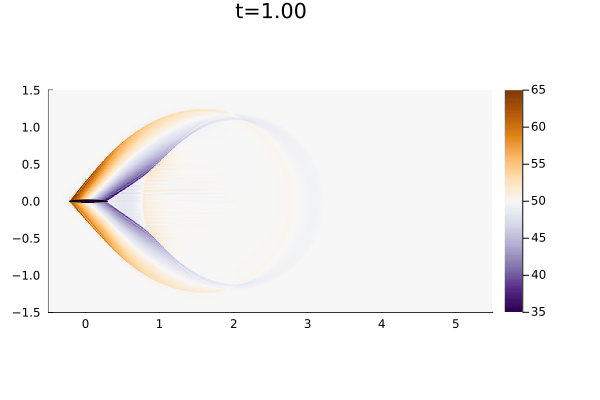}
    \end{subfigure}
    \vspace{7mm}
    
    \begin{subfigure}[b]{0.49\textwidth}
        \centering
        \includegraphics[trim={0cm 1cm 1.5cm 0cm}, clip, width=\textwidth]{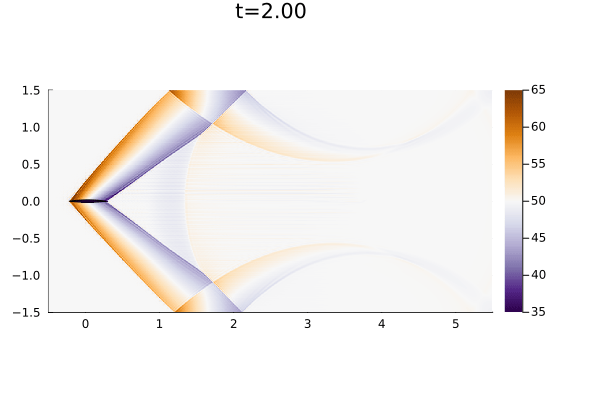}
    \end{subfigure}
    \begin{subfigure}[b]{0.49\textwidth}
        \centering
        \includegraphics[trim={0cm 1cm 1.5cm 0cm}, clip, width=\textwidth]{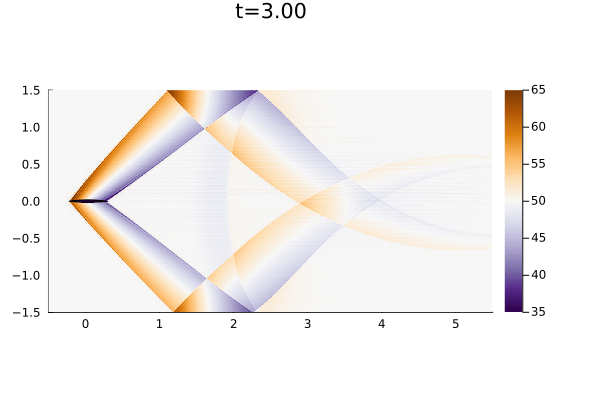}
    \end{subfigure}
    \vspace{7mm}

    \begin{subfigure}[b]{0.49\textwidth}
        \centering
        \includegraphics[trim={0cm 1cm 1.5cm 0cm}, clip, width=\textwidth]{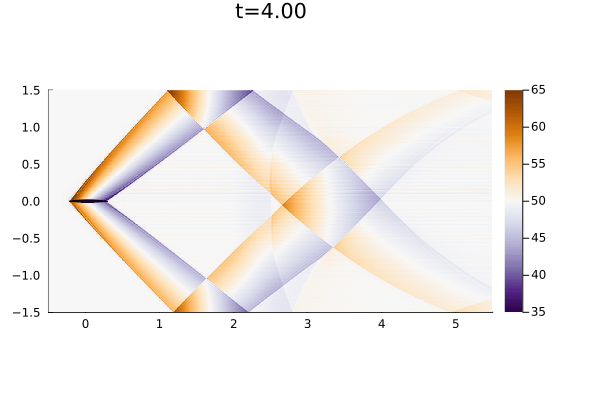}
    \end{subfigure}
    \begin{subfigure}[b]{0.49\textwidth}
        \centering
        \includegraphics[trim={0cm 1cm 1.5cm 0cm}, clip, width=\textwidth]{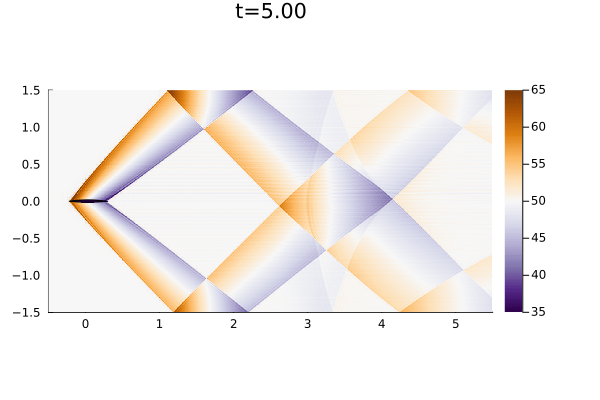}
    \end{subfigure}
    \caption{Snapshots of the biconvex airfoil simulation at various time stamps.}
    \label{fig:biconvex-snapshots}
\end{figure}

For the initial condition we use an impulsive start with the solution initialized to the freestream conditions. We take $N = 4$ and simulate the system up to time $t = 5$ with state redistribution using an initial time step of $\Delta t = 1 \times 10^{-8}$ (a typical timestep was on the order of $1 \times 10^{-4}$). During this time period, a bow shock immediately forms off the leading edge of the airfoil before reflecting off the top and bottom walls as shown in Figure \ref{fig:biconvex-snapshots}. Our method is able to capture these shock waves and shock-shock interactions without loss of stability.

\pagebreak
\subsection{Dam Break}
\begin{figure}[H]
    \centering
    \includegraphics[width=0.9\linewidth]{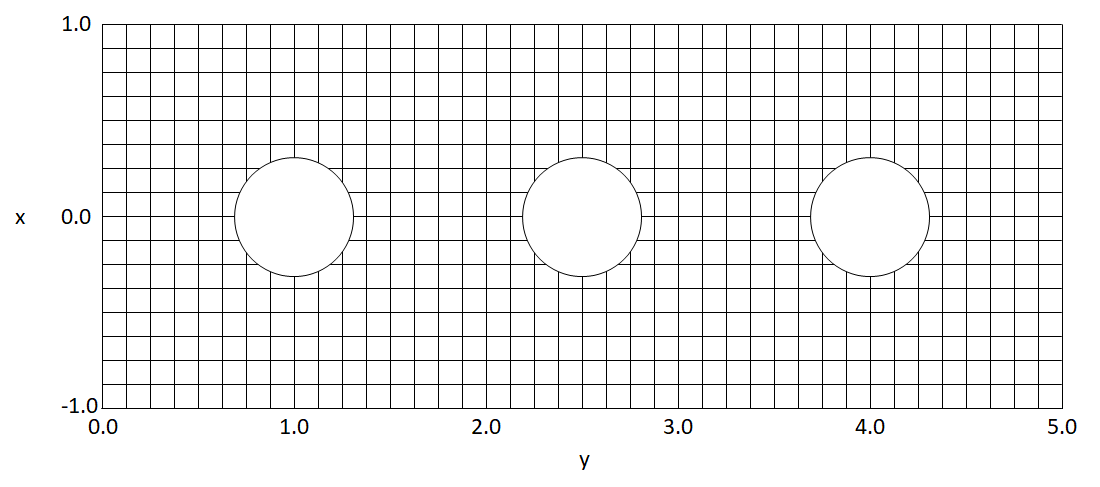}
    \caption{Mesh for the dam break experiment.}
    \label{fig:dam-break-mesh}
\end{figure}

For our last experiment we consider the shallow water equations on a domain with three embedded objects. We take the Cartesian embedding domain to be $\CartDomain = [-1,1] \times [0,5]$ over which we impose a $16 \times 40$ background Cartesian mesh. Figure \ref{fig:dam-break-mesh} shows the mesh for this problem. The three embedded circles share the same radius $R = 0.31$ and are centered at the points $(x,y) = (0,1), (0,2.5), (0.4)$. We again use a zero velocity initial condition and take the initial water height to be the discontinuous function
\begin{equation}
    h(x,y,t=0) = \left\{ \begin{matrix}
        3, &\quad y \geq 4.5 \\
        2, & \quad y < 4.5
    \end{matrix} \right. .
\end{equation}

\begin{figure}[H]
    \centering
    \begin{subfigure}[b]{0.45\textwidth}
        \centering
        \includegraphics[trim={4cm 1cm 3.5cm 2.5cm}, clip, width=\textwidth]{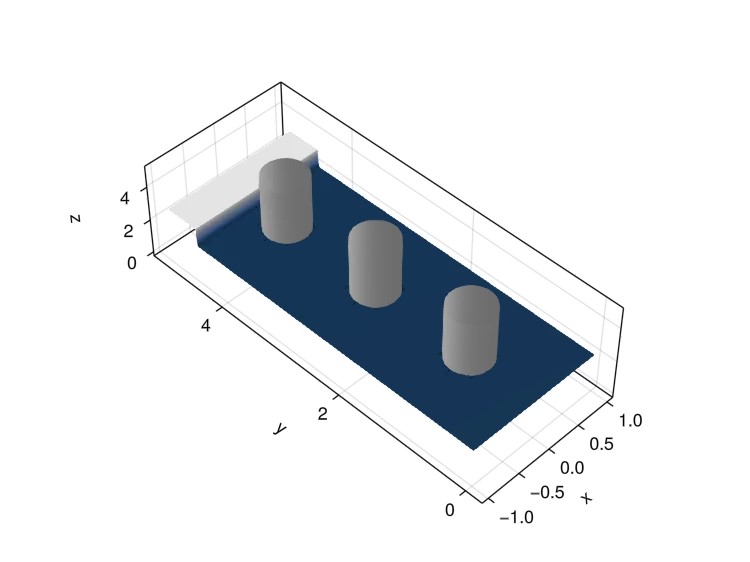}
        \caption{$t = 0.00$ }
    \end{subfigure}
    \begin{subfigure}[b]{0.45\textwidth}
        \centering
        \includegraphics[trim={4cm 1cm 3.5cm 2.5cm}, clip, width=\textwidth]{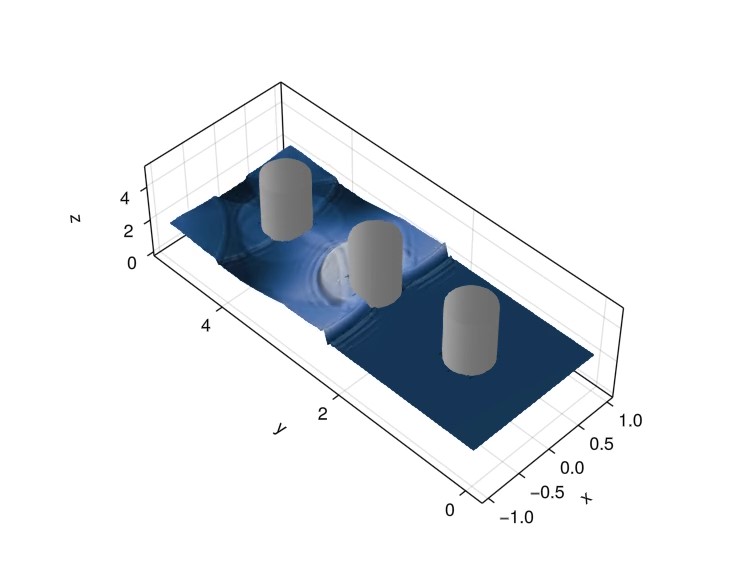}
        \caption{$t = 1.34$ }
    \end{subfigure}
    \begin{subfigure}[b]{0.45\textwidth}
        \centering
        \includegraphics[trim={4cm 1cm 3.5cm 2.5cm}, clip, width=\textwidth]{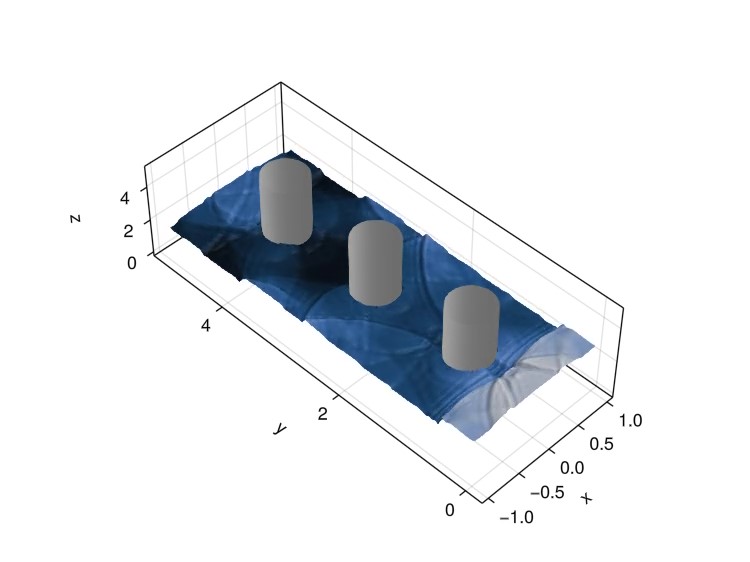}
        \caption{$t = 3.19$ }
    \end{subfigure}
    \begin{subfigure}[b]{0.45\textwidth}
        \centering
        \includegraphics[trim={4cm 1cm 3.5cm 2.5cm}, clip, width=\textwidth]{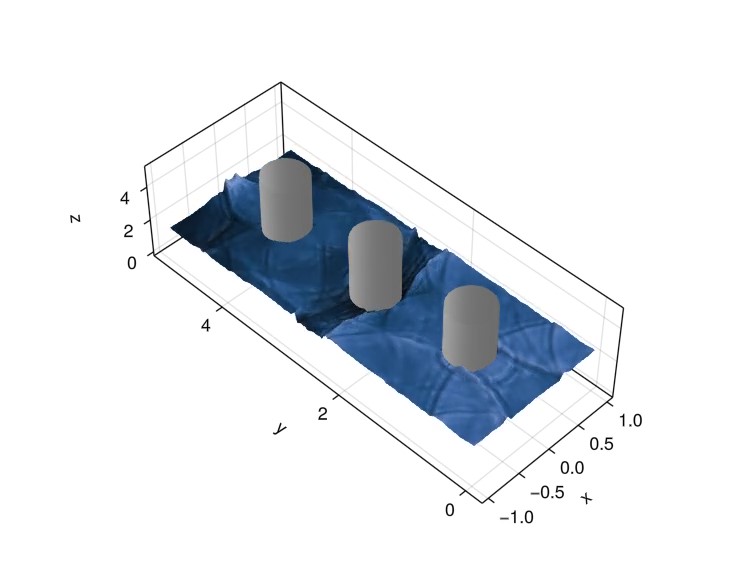}
        \caption{$t = 4.28$ }
    \end{subfigure}
    \begin{subfigure}[b]{0.45\textwidth}
        \centering
        \includegraphics[trim={4cm 1cm 3.5cm 2.5cm}, clip, width=\textwidth]{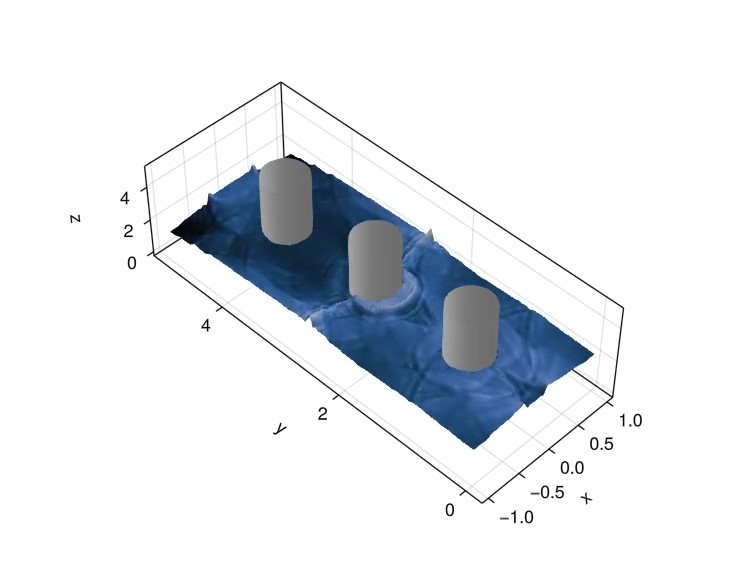}
        \caption{$t = 4.65$ }
    \end{subfigure}
    \begin{subfigure}[b]{0.45\textwidth}
        \centering
        \includegraphics[trim={4cm 1cm 3.5cm 2.5cm}, clip, width=\textwidth]{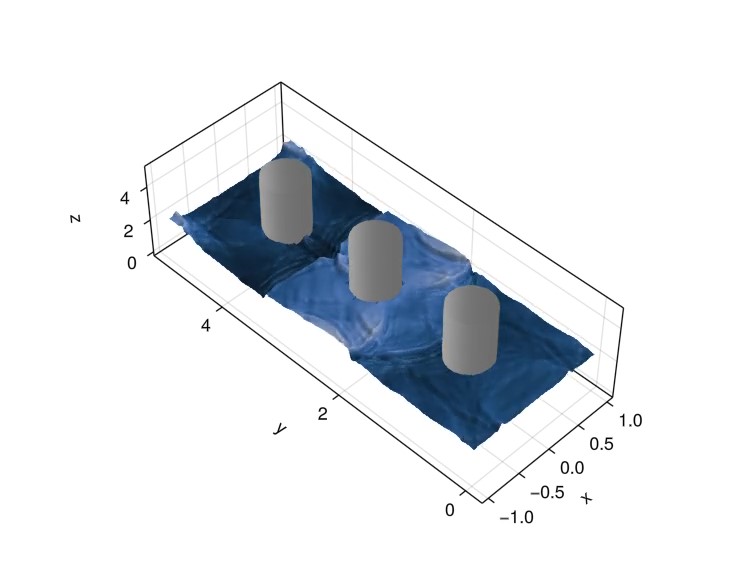}
        \caption{$t = 4.91$ }
    \end{subfigure}
    
    \caption{Snapshots of the dam break simulation at various times.}
    \label{fig:dam-break-snapshots}
\end{figure}

We impose reflective wall boundary conditions on all boundaries and take the solution degree as $N=4$. We simulate the system up to an end time of $t=5$ using an initial time step of $\Delta t = 1 \times 10^{-4}$ with state redistribution. This end time is sufficient for the front to propagate to and reflect off of the far wall. Figure \ref{fig:dam-break-snapshots} shows snapshots of the solution at various points in time.

\section{Conclusions}\label{sec:conclusions} 
We have demonstrated, to the authors' knowledge, the first high-order accurate entropy stable DG method on cut meshes which can be adapted to any hyperbolic conservation law taking the form in Equation \eqref{eq:hyperbolicConsLaw}. 

The presented method allows the robustness of entropy-stable schemes to be combined with the superior approximation power of high-order DG and efficiency of cut meshes. We combine these features via the skew-hybridized SBP formulation of Chan \cite{chan-SkewSymmES}, which allows the SBP property to be enforced on hybrid meshes for sufficiently accurate quadrature. 

To generate such quadrature on cut meshes we triangulate each cut element in the mesh and create an initial, many-point composite quadrature rule as the sum of quadrature rules on each subtriangle. These initial many-point rules are then pruned via Carath\'eodory pruning using the QR factorization based approach of van den Bos \cite{vandenbos-caratheodory}. We have numerically verified the entropy status and expected order of accuracy our method and demonstrated its accuracy and robustness on a number of problems. 

While we have not theoretically explored the implications of state redistribution on entropy stability, we have shown that state redistribution violates entropy conservation (as expected) but that, in the very least, it does not always violate entropy stability. In future work we would like to further explore the impact of state redistribution on entropy stable schemes to increase the robustness of time integration on cut meshes.

\section{Acknowledgments}
The authors gratefully acknowledge support by the National Science Foundation under the awards NSF GRFP-1842494, DMS-2231482, and DMS-1943186 and the Oden Institute for Computational Engineering and Science. The authors would also like to thank Akil Narayan for helpful conversations on Carath\'eodory pruning.

\clearpage

\section{Appendix}\label{ap:pruning}

For a given cut element $\element^k$ and target polynomial degree $M$, let $M^* = \dim (\totalDegree^M(\element^k))$ and let $\{\basisp_i\}_{i=1}^{M^*}$ be a basis for $\totalDegree^M(\element^k)$. Given an volume quadrature rule $(\bw, \bX)$ with $m$ points exact for polynomials in $\totalDegree^M(\element^k)$, we can express the quadrature rule as the linear system:
\begin{equation}
    \Phi(\bX) \bw = \bm{b}
\end{equation}
where for a set of points $\bX = \{\bx_1, ..., \bx_m \} \subset \realR^\Nspace$
\begin{equation}
    \Phi\LRp{\bX} = \begin{bmatrix}
        \basisp_1\LRp{\bx_1} & \cdots & \basisp_1\LRp{\bx_m} \\
        \vdots & & \vdots \\
        \basisp_{M^*}\LRp{\bx_1} & \cdots & \basisp_{M^*}\LRp{\bx_m}
    \end{bmatrix}, \quad \bm{b} = \begin{bmatrix}
        \int_{\element^k} \basisp_1 \\
        \vdots \\
        \int_{\element^k} \basisp_{M^*}
    \end{bmatrix}.
\end{equation}
Notice that $\Phi(\bX_q) = \bV_q^T$ when the same polynomial basis is used for both matrices (in practice the bases are different). Let $\Delta \bw$ be a null vector of $\Phi(\bX_q)$. By the definition of a null vector, for any scalar $\alpha$ we have that
\begin{equation}
    \Phi(\bX\LRp{\bw - \alpha\Delta\bw} = \Phi(\bX)\bw = \bm{b}.\label{eq:null-vector}
\end{equation}
Equation \eqref{eq:null-vector} is at the heart of quadrature pruning and/or modification techniques. It testifies that null vectors can be used to modify quadrature weights without impacting the exactness of the quadrature rule. In our case, we compute the null vector using a QR factorization of $\Phi(\bX)$. Once an arbitrary null vector of $\Phi\LRp{\bX}$ is known, we can choose the scalar $\alpha$ such that it zeros one or more quadrature weights while maintaining non-negativity in the other weights. We take $\alpha_0$ to be
\begin{equation}
    \alpha = \min \{\alpha^-, \alpha^+\} \label{eq:alpha-final}
\end{equation}
where
\begin{align}
    \alpha^- &= \max_{\Delta \bw_i < 0} \frac{\bw_i}{\Delta \bw_i}, \\[8pt]
    \alpha^+ &= \max_{\Delta \bw_i > 0} \frac{\bw_i}{\Delta \bw_i}.\label{eq:alpha+}
\end{align}
This definition for $\alpha$ guarantees that
\begin{equation}
    \bw_i - \alpha_i \Delta \bw_i \geq 0, \quad i = 1, ..., m,
\end{equation}
i.e., the quadrature weights all remain non-negative, with the entries used for $\alpha$ being zeroed. This process can be used to iteratively prune a quadrature rule one (or more) weights at a time. Let $(\bw^{(0)}, \bX^{(0)})$ be the initial $m$ point quadrature rule with non-negative quadrature weights. For a given stage of pruning let $\Delta \bw^{(i)}$ be a null vector of $\Phi(\bX^{(i)})$ and $\alpha_i$ defined as shown in \eqref{eq:alpha-final}-\eqref{eq:alpha+} using $\bw^{(i)}$ and $\Delta \bw^{(i)}$. Let $J = \{j_1, ..., j_{m_i}\}$ be the set of indices for all weights not zeroed by $\alpha_i \Delta \bw^{(i)}$. The updated quadrature rule is then given by
\begin{align}
    \bX^{(i+1)} &= \{ \bx_j \in \bX^{(i)} ~:~ j \in J\}\label{eq:pruning-step-X} \\[8pt]
    \bw^{(i+1)}_k &= \bw^{(i)}_{j_k}, \quad k = 1, ..., |J|. \label{eq:pruning-step-w}
\end{align}
The iteration is continued until a final quadrature rule with $M^* + 1$ points, the upper bound given by Carath\'eodory's theorem, or less is achieved.

The main cost of pruning comes from the QR-factorization, as it must be repeated at every iteration. It may be possible to greatly reduce the cost of this step via computing a structurally orthogonal basis for the null space of the $\Phi(\bX^{(0)})$ and thereby reducing the cost to a single QR factorizations and the orthongonalization routine. However, constructing such a routine is quite difficult (e.g., it can require choosing which weights should be pruned) and so for this proof of concept work, which only features static cut boundaries, we do not explore optimizing the pruning process.

However, since our initial quadrature rule is based on a subtriangulation, it is often the case that duplicate and/or near-duplicate quadrature points are present in the initial rules; such pairs of duplicate points arise on neighboring subtriangles on their shared face. When such duplicate points are present, it is very likely that multiple quadrature weights will be zeroed at once. In the case of near-duplicate points quadrature weights can also become arbitrarily small. As previously mentioned, non-negative quadrature weights are permissible so long as enough strictly positive weights remain for the mass matrix to be positive definite. In particular, the mass matrix will be positive definite if the quadrature points corresponding to strictly positive weights are sufficient to define a basis for $\totalDegree^N(\element^k)$. This condition is typically easily met: the quadrature rules are designed for polynomials of degree $2N-1$, and thus $n_q$ is typically much larger than $\dim(\totalDegree^N(\element^k))$ even after zeroed entries are discarded.

\bibliographystyle{elsarticle-num}
\bibliography{bibliography.bib}

\begin{thebibliography}{10}
\expandafter\ifx\csname url\endcsname\relax
  \def\url#1{\texttt{#1}}\fi
\expandafter\ifx\csname urlprefix\endcsname\relax\def\urlprefix{URL }\fi
\expandafter\ifx\csname href\endcsname\relax
  \def\href#1#2{#2} \def\path#1{#1}\fi

\bibitem{chan-SkewSymmES}
J.~Chan, Skew-symmetric entropy stable modal discontinuous {G}alerkin
  formulations, Journal of Scientific Computing 81~(1) (2019) 459--485.
\newblock \href {https://doi.org/10.1007/s10915-019-01026-w}
  {\path{doi:10.1007/s10915-019-01026-w}}.

\bibitem{dafermos-HyperbolicConservationLawsBook}
C.~M. Dafermos, Hyperbolic Conservation Laws in Continuum Physics, A Series of
  Comprehensive Studies in Mathematics, Springer-Verlag Berlin Heidelberg,
  2009.
\newblock \href {https://doi.org/10.1007/978-3-642-04048-1}
  {\path{doi:10.1007/978-3-642-04048-1}}.

\bibitem{gustafsson-L2-stability}
B.~Gustafsson, H.-O. Kreiss, J.~Oliger, Time dependent problems and difference
  methods, Vol.~24, John Wiley \& Sons, 1995.

\bibitem{wintermeyer-SWE-flux}
N.~Wintermeyer, A.~R. Winters, G.~J. Gassner, D.~A. Kopriva,
  \href{https://www.sciencedirect.com/science/article/pii/S0021999117302310}{An
  entropy stable nodal discontinuous {G}alerkin method for the two dimensional
  shallow water equations on unstructured curvilinear meshes with discontinuous
  bathymetry}, Journal of Computational Physics 340 (2017) 200--242.
\newblock \href {https://doi.org/10.1016/j.jcp.2017.03.036}
  {\path{doi:10.1016/j.jcp.2017.03.036}}.
\newline\urlprefix\url{https://www.sciencedirect.com/science/article/pii/S0021999117302310}

\bibitem{harten-ES-gudonov}
A.~Harten, P.~D. Lax, B.~v. Leer, \href{https://doi.org/10.1137/1025002}{On
  upstream differencing and {G}odunov-type schemes for hyperbolic conservation
  laws}, SIAM Review 25~(1) (1983) 35--61.
\newblock \href {https://doi.org/10.1137/1025002} {\path{doi:10.1137/1025002}}.
\newline\urlprefix\url{https://doi.org/10.1137/1025002}

\bibitem{chen-ESFluxes}
T.~Chen, C.-W. Shu, Entropy stable high order discontinuous {G}alerkin methods
  with suitable quadrature rules for hyperbolic conservation laws, Journal of
  Computational Physics 345 (2017) 427--461.
\newblock \href {https://doi.org/10.1016/j.jcp.2017.05.025}
  {\path{doi:10.1016/j.jcp.2017.05.025}}.

\bibitem{godlewski-hyperbolicMethodsBook}
E.~Godlewski, P.-A. Raviart, Numerical approximation of hyperbolic systems of
  conservation laws, Vol. 118, Springer Science \& Business Media, 2013.

\bibitem{hughes-NavierStokesFlux}
T.~Hughes, L.~Franca, M.~Mallet,
  \href{https://www.sciencedirect.com/science/article/pii/0045782586901271}{A
  new finite element formulation for computational fluid dynamics: I. symmetric
  forms of the compressible {E}uler and {N}avier-{S}tokes equations and the
  second law of thermodynamics}, Computer Methods in Applied Mechanics and
  Engineering 54~(2) (1986) 223--234.
\newblock \href {https://doi.org/10.1016/0045-7825(86)90127-1}
  {\path{doi:10.1016/0045-7825(86)90127-1}}.
\newline\urlprefix\url{https://www.sciencedirect.com/science/article/pii/0045782586901271}

\bibitem{harten-ESChainrule}
A.~Harten, On the symmetric form of systems of conservation laws with entropy,
  Journal of Computational Physics 49~(1) (1983) 151--164.
\newblock \href {https://doi.org/10.1016/0021-9991(83)90118-3}
  {\path{doi:10.1016/0021-9991(83)90118-3}}.

\bibitem{tadmor-fluxes}
E.~Tadmor, The numerical viscosity of entropy stable schemes for systems of
  conservation laws. {I}, Mathematics of Computation 49~(179) (1987) 91--103.

\bibitem{tadmor-EStheory}
E.~Tadmor, Entropy stability theory for difference approximations of nonlinear
  conservation laws and related time-dependent problems, Acta Numerica 12
  (2003) 451--512.

\bibitem{tadmor-ES-FV-FD-schemes}
E.~Tadmor,
  \href{https://www.sciencedirect.com/science/article/pii/S1570865916300151}{Chapter
  18 - {E}ntropy stable schemes}, in: R.~Abgrall, C.-W. Shu (Eds.), Handbook of
  Numerical Methods for Hyperbolic Problems, Vol.~17 of Handbook of Numerical
  Analysis, Elsevier, 2016, pp. 467--493.
\newblock \href {https://doi.org/10.1016/bs.hna.2016.09.006}
  {\path{doi:10.1016/bs.hna.2016.09.006}}.
\newline\urlprefix\url{https://www.sciencedirect.com/science/article/pii/S1570865916300151}

\bibitem{ray-ES-FV}
D.~Ray, P.~Chandrashekar, U.~S. Fjordholm, S.~Mishra, Entropy stable scheme on
  two-dimensional unstructured grids for {E}uler equations, Communications in
  Computational Physics 19~(5) (2016) 1111–1140.
\newblock \href {https://doi.org/10.4208/cicp.scpde14.43s}
  {\path{doi:10.4208/cicp.scpde14.43s}}.

\bibitem{fjordholm-TeCNO}
U.~S. Fjordholm, S.~Mishra, E.~Tadmor,
  \href{https://doi.org/10.1137/110836961}{Arbitrarily high-order accurate
  entropy stable essentially nonoscillatory schemes for systems of conservation
  laws}, SIAM Journal on Numerical Analysis 50~(2) (2012) 544--573.
\newblock \href {https://doi.org/10.1137/110836961}
  {\path{doi:10.1137/110836961}}.
\newline\urlprefix\url{https://doi.org/10.1137/110836961}

\bibitem{chandrashekar-KEP-ES-FV}
P.~Chandrashekar, Kinetic energy preserving and entropy stable finite volume
  schemes for compressible {E}uler and {N}avier-{S}tokes equations,
  Communications in Computational Physics 14~(5) (2013) 1252–1286.
\newblock \href {https://doi.org/10.4208/cicp.170712.010313a}
  {\path{doi:10.4208/cicp.170712.010313a}}.

\bibitem{harten-HLL-solver}
A.~Harten, P.~D. Lax, B.~v. Leer, \href{10.1137/1025002}{On upstream
  differencing and {G}odunov-type schemes for hyperbolic conservation laws},
  SIAM Review 25~(1) (1983) 35--61.
\newblock \href {https://doi.org/10.1137/1025002} {\path{doi:10.1137/1025002}}.
\newline\urlprefix\url{10.1137/1025002}

\bibitem{tan-laxWendroff}
S.~Tan, C.-W. Shu,
  \href{https://www.sciencedirect.com/science/article/pii/S0021999110003979}{Inverse
  {L}ax-{W}endroff procedure for numerical boundary conditions of conservation
  laws}, Journal of Computational Physics 229~(21) (2010) 8144--8166.
\newblock \href {https://doi.org/10.1016/j.jcp.2010.07.014}
  {\path{doi:10.1016/j.jcp.2010.07.014}}.
\newline\urlprefix\url{https://www.sciencedirect.com/science/article/pii/S0021999110003979}

\bibitem{guermond-lambda-euler}
J.-L. Guermond, B.~Popov,
  \href{https://www.sciencedirect.com/science/article/pii/S0021999116301991}{Fast
  estimation from above of the maximum wave speed in the {R}iemann problem for
  the {E}uler equations}, Journal of Computational Physics 321 (2016) 908--926.
\newblock \href {https://doi.org/10.1016/j.jcp.2016.05.054}
  {\path{doi:10.1016/j.jcp.2016.05.054}}.
\newline\urlprefix\url{https://www.sciencedirect.com/science/article/pii/S0021999116301991}

\bibitem{ainsworth-highOrder-dispersionError}
M.~Ainsworth, Dispersive and dissipative behaviour of high order discontinuous
  {G}alerkin finite element methods, Journal of Computational Physics 198~(1)
  (2004) 106--130.
\newblock \href {https://doi.org/10.1016/j.jcp.2004.01.004}
  {\path{doi:10.1016/j.jcp.2004.01.004}}.

\bibitem{wang-highVsLowOrder-HOinstability}
Z.~Wang, K.~Fidkowski, R.~Abgrall, F.~Bassi, D.~Caraeni, A.~Cary, H.~Deconinck,
  R.~Hartmann, K.~Hillewaert, H.~Huynh, N.~Kroll, G.~May, P.-O. Persson, B.~van
  Leer, M.~Visbal, High-order {CFD} methods: current status and perspective,
  International Journal for Numerical Methods in Fluids 72~(8) (2013) 811--845.
\newblock \href {https://doi.org/10.1002/fld.3767}
  {\path{doi:10.1002/fld.3767}}.

\bibitem{visbal-highOrder-unsteadyFlows}
M.~R. Visbal, D.~V. Gaitonde, High-order-accurate methods for complex unsteady
  subsonic flows, AIAA Journal 37~(10) (1999) 1231--1239.
\newblock \href {https://doi.org/10.2514/2.591} {\path{doi:10.2514/2.591}}.

\bibitem{tominec-blastWave}
I.~Tominec, M.~Nazarov,
  \href{https://doi.org/10.1007/s10915-022-02055-8}{Residual viscosity
  stabilized {RBF-FD} methods for solving nonlinear conservation laws}, Journal
  of Scientific Computing 94~(1) (2022) 14.
\newblock \href {https://doi.org/10.1007/s10915-022-02055-8}
  {\path{doi:10.1007/s10915-022-02055-8}}.
\newline\urlprefix\url{https://doi.org/10.1007/s10915-022-02055-8}

\bibitem{nazarov-blastWave-original}
M.~Nazarov, A.~Larcher,
  \href{https://www.sciencedirect.com/science/article/pii/S0045782516317704}{Numerical
  investigation of a viscous regularization of the {E}uler equations by entropy
  viscosity}, Computer Methods in Applied Mechanics and Engineering 317 (2017)
  128--152.
\newblock \href {https://doi.org/10.1016/j.cma.2016.12.010}
  {\path{doi:10.1016/j.cma.2016.12.010}}.
\newline\urlprefix\url{https://www.sciencedirect.com/science/article/pii/S0045782516317704}

\bibitem{zhang-maximumPrinciple}
X.~Zhang, C.-W. Shu,
  \href{https://www.sciencedirect.com/science/article/pii/S0021999109007165}{On
  maximum-principle-satisfying high order schemes for scalar conservation
  laws}, Journal of Computational Physics 229~(9) (2010) 3091--3120.
\newblock \href {https://doi.org/10.1016/j.jcp.2009.12.030}
  {\path{doi:10.1016/j.jcp.2009.12.030}}.
\newline\urlprefix\url{https://www.sciencedirect.com/science/article/pii/S0021999109007165}

\bibitem{zhang-positivityPreservation}
X.~Zhang, C.-W. Shu,
  \href{https://www.sciencedirect.com/science/article/pii/S0021999110004535}{On
  positivity-preserving high order discontinuous {G}alerkin schemes for
  compressible {E}uler equations on rectangular meshes}, Journal of
  Computational Physics 229~(23) (2010) 8918--8934.
\newblock \href {https://doi.org/10.1016/j.jcp.2010.08.016}
  {\path{doi:10.1016/j.jcp.2010.08.016}}.
\newline\urlprefix\url{https://www.sciencedirect.com/science/article/pii/S0021999110004535}

\bibitem{chan-hybridizedSBP}
J.~Chan, On discretely entropy conservative and entropy stable discontinuous
  {G}alerkin methods, Journal of Computational Physics 362 (2018) 346--374.
\newblock \href {https://doi.org/10.1016/j.jcp.2018.02.033}
  {\path{doi:10.1016/j.jcp.2018.02.033}}.

\bibitem{fisher-ES-FD}
T.~C. Fisher, M.~H. Carpenter,
  \href{https://www.sciencedirect.com/science/article/pii/S0021999113004385}{High-order
  entropy stable finite difference schemes for nonlinear conservation laws:
  Finite domains}, Journal of Computational Physics 252 (2013) 518--557.
\newblock \href {https://doi.org/10.1016/j.jcp.2013.06.014}
  {\path{doi:10.1016/j.jcp.2013.06.014}}.
\newline\urlprefix\url{https://www.sciencedirect.com/science/article/pii/S0021999113004385}

\bibitem{carpenter-ES-spectralCollocation}
M.~H. Carpenter, T.~C. Fisher, E.~J. Nielsen, S.~H. Frankel,
  \href{https://doi.org/10.1137/130932193}{Entropy stable spectral collocation
  schemes for the {N}avier--{S}tokes equations: Discontinuous interfaces}, SIAM
  Journal on Scientific Computing 36~(5) (2014) B835--B867.
\newblock \href {https://doi.org/10.1137/130932193}
  {\path{doi:10.1137/130932193}}.
\newline\urlprefix\url{https://doi.org/10.1137/130932193}

\bibitem{gassner-splitFormDG}
G.~J. Gassner, A.~R. Winters, D.~A. Kopriva,
  \href{https://www.sciencedirect.com/science/article/pii/S0021999116304259}{Split
  form nodal discontinuous {G}alerkin schemes with summation-by-parts property
  for the compressible {E}uler equations}, Journal of Computational Physics 327
  (2016) 39--66.
\newblock \href {https://doi.org/10.1016/j.jcp.2016.09.013}
  {\path{doi:10.1016/j.jcp.2016.09.013}}.
\newline\urlprefix\url{https://www.sciencedirect.com/science/article/pii/S0021999116304259}

\bibitem{gassner-BR1}
G.~J. Gassner, A.~R. Winters, F.~J. Hindenlang, D.~A. Kopriva,
  \href{https://doi.org/10.1007/s10915-018-0702-1}{The {BR1} scheme is stable
  for the compressible {N}avier--{S}tokes equations}, Journal of Scientific
  Computing 77~(1) (2018) 154--200.
\newblock \href {https://doi.org/10.1007/s10915-018-0702-1}
  {\path{doi:10.1007/s10915-018-0702-1}}.
\newline\urlprefix\url{https://doi.org/10.1007/s10915-018-0702-1}

\bibitem{crean-SBP-generalCurvedElements}
J.~Crean, J.~E. Hicken, D.~C. {Del Rey Fernández}, D.~W. Zingg, M.~H.
  Carpenter,
  \href{https://www.sciencedirect.com/science/article/pii/S0021999117308999}{Entropy-stable
  summation-by-parts discretization of the {E}uler equations on general curved
  elements}, Journal of Computational Physics 356 (2018) 410--438.
\newblock \href {https://doi.org/10.1016/j.jcp.2017.12.015}
  {\path{doi:10.1016/j.jcp.2017.12.015}}.
\newline\urlprefix\url{https://www.sciencedirect.com/science/article/pii/S0021999117308999}

\bibitem{reed-originalCutCell}
W.~H. Reed, T.~R. Hill, Triangular mesh methods for the neutron transport
  equation, Tech. Rep. {LA-UR-73-479}, Los Alamos Scientific Lab., N. Mex.(USA)
  (1973).

\bibitem{delReyFernandez-SBPframework}
D.~C. {Del Rey Fernández}, P.~D. Boom, D.~W. Zingg,
  \href{https://www.sciencedirect.com/science/article/pii/S002199911400076X}{A
  generalized framework for nodal first derivative summation-by-parts
  operators}, Journal of Computational Physics 266 (2014) 214--239.
\newblock \href {https://doi.org/10.1016/j.jcp.2014.01.038}
  {\path{doi:10.1016/j.jcp.2014.01.038}}.
\newline\urlprefix\url{https://www.sciencedirect.com/science/article/pii/S002199911400076X}

\bibitem{chen-ES-review}
T.~Chen, C.-W. Shu, Review of entropy stable discontinuous {G}alerkin methods
  for systems of conservation laws on unstructured simplex meshes, CSIAM
  Transactions on Applied Mathematics 1~(1) (2020) 1--52.

\bibitem{chan-collocation}
J.~Chan, D.~C. Del Rey~Fern\'{a}ndez, M.~H. Carpenter,
  \href{https://doi.org/10.1137/18M1209234}{Efficient entropy stable {G}auss
  collocation methods}, SIAM Journal on Scientific Computing 41~(5) (2019)
  A2938--A2966.
\newblock \href {https://doi.org/10.1137/18M1209234}
  {\path{doi:10.1137/18M1209234}}.
\newline\urlprefix\url{https://doi.org/10.1137/18M1209234}

\bibitem{stavrev-knownShapeCutCells}
A.~Stavrev, L.~H. Nguyen, R.~Shen, V.~Varduhn, M.~Behr, S.~Elgeti,
  D.~Schillinger,
  \href{https://www.sciencedirect.com/science/article/pii/S004578251630411X}{Geometrically
  accurate, efficient, and flexible quadrature techniques for the tetrahedral
  finite cell method}, Computer Methods in Applied Mechanics and Engineering
  310 (2016) 646--673.
\newblock \href {https://doi.org/10.1016/j.cma.2016.07.041}
  {\path{doi:10.1016/j.cma.2016.07.041}}.
\newline\urlprefix\url{https://www.sciencedirect.com/science/article/pii/S004578251630411X}

\bibitem{davis-momentFitting}
P.~J. Davis, A construction of nonnegative approximate quadratures, Mathematics
  of Computation 21~(100) (1967) 578--582.
\newblock \href {https://doi.org/10.2307/2005001} {\path{doi:10.2307/2005001}}.

\bibitem{garhoum-momentFitting1}
W.~Garhuom, A.~D{\"u}ster,
  \href{https://doi.org/10.1007/s00466-022-02203-9}{Non-negative moment fitting
  quadrature for cut finite elements and cells undergoing large deformations},
  Computational Mechanics 70~(5) (2022) 1059--1081.
\newblock \href {https://doi.org/10.1007/s00466-022-02203-9}
  {\path{doi:10.1007/s00466-022-02203-9}}.
\newline\urlprefix\url{https://doi.org/10.1007/s00466-022-02203-9}

\bibitem{bui-momentFitting2}
H.-G. Bui, D.~Schillinger, G.~Meschke,
  \href{https://www.sciencedirect.com/science/article/pii/S0045782520302346}{Efficient
  cut-cell quadrature based on moment fitting for materially nonlinear
  analysis}, Computer Methods in Applied Mechanics and Engineering 366 (2020)
  113050.
\newblock \href {https://doi.org/10.1016/j.cma.2020.113050}
  {\path{doi:10.1016/j.cma.2020.113050}}.
\newline\urlprefix\url{https://www.sciencedirect.com/science/article/pii/S0045782520302346}

\bibitem{legrain-momentFitting}
G.~Legrain,
  \href{https://www.sciencedirect.com/science/article/pii/S0898122121002820}{Non-negative
  moment fitting quadrature rules for fictitious domain methods}, Computers \&
  Mathematics with Applications 99 (2021) 270--291.
\newblock \href {https://doi.org/10.1016/j.camwa.2021.07.019}
  {\path{doi:10.1016/j.camwa.2021.07.019}}.
\newline\urlprefix\url{https://www.sciencedirect.com/science/article/pii/S0898122121002820}

\bibitem{piazzon-opt-pruning}
F.~Piazzon, A.~Sommariva, M.~Vianello, et~al., {C}aratheodory-{T}chakaloff
  least squares, in: International Conference on Sampling Theory and
  Applications (SampTA), 2017, pp. 672--676.

\bibitem{saye-Quadrature}
R.~I. Saye, High-order quadrature methods for implicitly defined surfaces and
  volumes in hyperrectangles, SIAM Journal on Scientific Computing 37~(2)
  (2015) A993--A1019.
\newblock \href {https://doi.org/10.1137/140966290}
  {\path{doi:10.1137/140966290}}.

\bibitem{saye-QuadraturePolynomials}
R.~I. Saye, High-order quadrature on multi-component domains implicitly defined
  by multivariate polynomials, Journal of Computational Physics 448 (2022)
  110720.
\newblock \href {https://doi.org/10.1016/j.jcp.2021.110720}
  {\path{doi:10.1016/j.jcp.2021.110720}}.

\bibitem{kaur-DGDiff}
S.~Kaur, G.~Yan, J.~E. Hicken, High-order cut-cell discontinuous {G}alerkin
  difference discretization, AIAA Journal 61~(10) (2023) 4220--4229.
\newblock \href {https://doi.org/10.2514/1.J062990}
  {\path{doi:10.2514/1.J062990}}.

\bibitem{taylor-L2-cutDG}
C.~G. Taylor, L.~C. Wilcox, J.~Chan,
  \href{https://www.sciencedirect.com/science/article/pii/S0021999124007769}{An
  energy stable high-order cut cell discontinuous {G}alerkin method with state
  redistribution for wave propagation}, Journal of Computational Physics 521
  (2025) 113528.
\newblock \href {https://doi.org/10.1016/j.jcp.2024.113528}
  {\path{doi:10.1016/j.jcp.2024.113528}}.
\newline\urlprefix\url{https://www.sciencedirect.com/science/article/pii/S0021999124007769}

\bibitem{kudela-quadtreeQuadr}
L.~Kudela, N.~Zander, T.~Bog, S.~Kollmannsberger, E.~Rank, Efficient and
  accurate numerical quadrature for immersed boundary methods, Advanced
  modeling and simulation in engineering sciences 2 (2015) 1--22.

\bibitem{burman-cutFEM1}
E.~Burman, S.~Claus, P.~Hansbo, M.~G. Larson, A.~Massing,
  \href{https://onlinelibrary.wiley.com/doi/abs/10.1002/nme.4823}{{CutFEM}:
  Discretizing geometry and partial differential equations}, International
  Journal for Numerical Methods in Engineering 104~(7) (2015) 472--501.
\newblock \href {https://doi.org/10.1002/nme.4823}
  {\path{doi:10.1002/nme.4823}}.
\newline\urlprefix\url{https://onlinelibrary.wiley.com/doi/abs/10.1002/nme.4823}

\bibitem{burman-cutFEM2}
E.~Burman, P.~Hansbo, M.~G. Larson, {CutFEM} based on extended finite element
  spaces, Numerische Mathematik 152~(2) (2022) 331--369.

\bibitem{hernandez-quadrature}
J.~Hernández, M.~Caicedo, A.~Ferrer, Dimensional hyper-reduction of nonlinear
  finite element models via empirical cubature, Computer Methods in Applied
  Mechanics and Engineering 313 (2017) 687--722.
\newblock \href {https://doi.org/10.1016/j.cma.2016.10.022}
  {\path{doi:10.1016/j.cma.2016.10.022}}.

\bibitem{slobodkins-quadrature-nodeElimination}
A.~Slobodkins, J.~Tausch,
  \href{https://www.sciencedirect.com/science/article/pii/S0898122123002468}{A
  node elimination algorithm for cubature of high-dimensional polytopes},
  Computers \& Mathematics with Applications 144 (2023) 229--236.
\newblock \href {https://doi.org/10.1016/j.camwa.2023.06.001}
  {\path{doi:10.1016/j.camwa.2023.06.001}}.
\newline\urlprefix\url{https://www.sciencedirect.com/science/article/pii/S0898122123002468}

\bibitem{vandenbos-caratheodory}
L.~{van den Bos}, B.~Sanderse, W.~Bierbooms, G.~{van Bussel},
  \href{https://epubs.siam.org/doi/abs/10.1137/18M1213373}{Generating nested
  quadrature rules with positive weights based on arbitrary sample sets},
  SIAM/ASA Journal on Uncertainty Quantification 8~(1) (2020) 139--169.
\newblock \href {https://doi.org/10.1137/18M1213373}
  {\path{doi:10.1137/18M1213373}}.
\newline\urlprefix\url{https://epubs.siam.org/doi/abs/10.1137/18M1213373}

\bibitem{wilson-caratheodory}
M.~W. Wilson,
  \href{https://www.ams.org/mcom/1969-23-106/S0025-5718-1969-0242374-1/}{A
  general algorithm for nonnegative quadrature formulas}, Mathematics of
  Computation 23~(106) (1969) 253--258.
\newblock \href {https://doi.org/10.1090/S0025-5718-1969-0242374-1}
  {\path{doi:10.1090/S0025-5718-1969-0242374-1}}.
\newline\urlprefix\url{https://www.ams.org/mcom/1969-23-106/S0025-5718-1969-0242374-1/}

\bibitem{caratheodory-thm}
C.~Carath{\'e}odory, \href{https://doi.org/10.1007/BF03014795}{{\"U}ber den
  variabilit{\"a}tsbereich der fourier'schen konstanten von positiven
  harmonischen funktionen}, Rendiconti del Circolo Matematico di Palermo
  (1884-1940) 32~(1) (1911) 193--217.
\newblock \href {https://doi.org/10.1007/BF03014795}
  {\path{doi:10.1007/BF03014795}}.
\newline\urlprefix\url{https://doi.org/10.1007/BF03014795}

\bibitem{vioreanu-quadrature}
B.~Vioreanu, V.~Rokhlin, \href{https://doi.org/10.1137/110860082}{Spectra of
  multiplication operators as a numerical tool}, SIAM Journal on Scientific
  Computing 36~(1) (2014) A267--A288.
\newblock \href {https://doi.org/10.1137/110860082}
  {\path{doi:10.1137/110860082}}.
\newline\urlprefix\url{https://doi.org/10.1137/110860082}

\bibitem{Julia-PathIntersections}
\href{https://github.com/cgt3/PathIntersections.jl}{{PathIntersections.jl}}.
\newline\urlprefix\url{https://github.com/cgt3/PathIntersections.jl}

\bibitem{Julia-StartUpDG}
\href{https://github.com/jlchan/StartUpDG.jl}{{StartUpDG.jl}}.
\newline\urlprefix\url{https://github.com/jlchan/StartUpDG.jl}

\bibitem{code-repo}
\href{https://github.com/cgt3/ES-CutDG}{Reproducability repository for
  simulation codes}.
\newline\urlprefix\url{https://github.com/cgt3/ES-CutDG}

\bibitem{berger-stateRedistr}
M.~Berger, A.~Giuliani, A state redistribution algorithm for finite volume
  schemes on cut cell meshes, Journal of Computational Physics 428 (2021)
  109820.
\newblock \href {https://doi.org/10.1016/j.jcp.2020.109820}
  {\path{doi:10.1016/j.jcp.2020.109820}}.

\bibitem{giuliani-DG-SRD}
A.~Giuliani, \href{https://doi.org/10.1137/21M1396277}{A two-dimensional
  stabilized discontinuous {G}alerkin method on curvilinear embedded boundary
  grids}, SIAM Journal on Scientific Computing 44~(1) (2022) A389--A415.
\newblock \href {https://doi.org/10.1137/21M1396277}
  {\path{doi:10.1137/21M1396277}}.
\newline\urlprefix\url{https://doi.org/10.1137/21M1396277}

\bibitem{giuliani-weightedSRD}
A.~Giuliani, A.~Almgren, J.~Bell, M.~Berger, M.~{Henry de Frahan},
  D.~Rangarajan, A weighted state redistribution algorithm for embedded
  boundary grids, Journal of Computational Physics 464 (2022) 111305.
\newblock \href {https://doi.org/10.1016/j.jcp.2022.111305}
  {\path{doi:10.1016/j.jcp.2022.111305}}.

\bibitem{berger-weightedSRD}
M.~Berger, A.~Giuliani, A new provably stable weighted state redistribution
  algorithm, SIAM Journal on Scientific Computing 46~(5) (2024) A2848--A2873.
\newblock \href {https://doi.org/10.1137/23M1597484}
  {\path{doi:10.1137/23M1597484}}.

\bibitem{sommariva-feketeQR}
A.~Sommariva, M.~Vianello, Computing approximate {F}ekete points by {QR}
  factorizations of {V}andermonde matrices, Computers \& Mathematics with
  Applications 57~(8) (2009) 1324--1336.
\newblock \href {https://doi.org/10.1016/j.camwa.2008.11.011}
  {\path{doi:10.1016/j.camwa.2008.11.011}}.

\bibitem{shi-shu-localConservation}
C.~Shi, C.-W. Shu,
  \href{https://www.sciencedirect.com/science/article/pii/S004579301730230X}{On
  local conservation of numerical methods for conservation laws}, Computers \&
  Fluids 169 (2018) 3--9.
\newblock \href {https://doi.org/10.1016/j.compfluid.2017.06.018}
  {\path{doi:10.1016/j.compfluid.2017.06.018}}.
\newline\urlprefix\url{https://www.sciencedirect.com/science/article/pii/S004579301730230X}

\bibitem{davis-wavespeed}
S.~F. Davis, \href{https://doi.org/10.1137/0909030}{Simplified second-order
  {G}odunov-type methods}, SIAM Journal on Scientific and Statistical Computing
  9~(3) (1988) 445--473.
\newblock \href {https://doi.org/10.1137/0909030} {\path{doi:10.1137/0909030}}.
\newline\urlprefix\url{https://doi.org/10.1137/0909030}

\bibitem{ranocha-EC-eulerFluxes}
H.~Ranocha, \href{https://doi.org/10.1007/s10915-017-0618-1}{Comparison of some
  entropy conservative numerical fluxes for the {E}uler equations}, Journal of
  Scientific Computing 76~(1) (2018) 216--242.
\newblock \href {https://doi.org/10.1007/s10915-017-0618-1}
  {\path{doi:10.1007/s10915-017-0618-1}}.
\newline\urlprefix\url{https://doi.org/10.1007/s10915-017-0618-1}

\bibitem{ranocha-EulerFlux}
H.~Ranocha, Entropy conserving and kinetic energy preserving numerical methods
  for the {E}uler equations using summation-by-parts operators, in: S.~J.
  Sherwin, D.~Moxey, J.~Peir{\'o}, P.~E. Vincent, C.~Schwab (Eds.), Spectral
  and High Order Methods for Partial Differential Equations ICOSAHOM 2018,
  Springer International Publishing, Cham, 2020, pp. 525--535.

\bibitem{tsitouras-Tsit5}
C.~Tsitouras, {Runge–Kutta} pairs of order 5(4) satisfying only the first
  column simplifying assumption, Computers \& Mathematics with Applications
  62~(2) (2011) 770--775.
\newblock \href {https://doi.org/10.1016/j.camwa.2011.06.002}
  {\path{doi:10.1016/j.camwa.2011.06.002}}.

\bibitem{Julia-OrdinaryDiff}
C.~Rackauckas, Q.~Nie, {D}ifferential{E}quations.jl--a performant and
  feature-rich ecosystem for solving differential equations in {J}ulia, Journal
  of Open Research Software 5~(1) (2017) 15--15.
\newblock \href {https://doi.org/10.5334/jors.151}
  {\path{doi:10.5334/jors.151}}.

\bibitem{svard-ES-eulerBCs}
M.~Sv{\"a}rd, H.~{\"O}zcan,
  \href{https://doi.org/10.1007/s10915-013-9727-7}{Entropy-stable schemes for
  the {E}uler equations with far-field and wall boundary conditions}, Journal
  of Scientific Computing 58~(1) (2014) 61--89.
\newblock \href {https://doi.org/10.1007/s10915-013-9727-7}
  {\path{doi:10.1007/s10915-013-9727-7}}.
\newline\urlprefix\url{https://doi.org/10.1007/s10915-013-9727-7}

\bibitem{chu-entropyWave}
B.-T. Chu, L.~S.~G. Kovásznay, Non-linear interactions in a viscous
  heat-conducting compressible gas, Journal of Fluid Mechanics 3~(5) (1958)
  494–514.
\newblock \href {https://doi.org/10.1017/S0022112058000148}
  {\path{doi:10.1017/S0022112058000148}}.

\end{thebibliography}

\end{document}